\documentclass[11pt]{article}
\usepackage{amsmath,amssymb,amsthm}
\usepackage{dsfont}
\usepackage{graphicx}
\usepackage{geometry}
\usepackage[british]{babel}
\usepackage[utf8]{inputenc}
\usepackage{xcolor}
\usepackage[linesnumbered,ruled,vlined,commentsnumbered]{algorithm2e}
\usepackage{cite}
\usepackage{url}
\usepackage{diagbox}
\usepackage{hyperref}

\def\vs{\lambda_S}
\def\ve{\lambda_E}
\def\vi{\lambda_I}
\def\va{\lambda_A}
\def\vr{\lambda_R}
\def\RR{\mathbb{R}}

\newtheorem{remark}{Remark}

\begin{document}
\title{Bi-fidelity stochastic collocation methods for epidemic transport models with uncertainties}

\author{Giulia Bertaglia\footnote{Department of Mathematics and Computer Science, University of Ferrara, Via Machiavelli 30, 44121 - Ferrara, Italy ({\tt giulia.bertaglia@unife.it}, {\tt lorenzo.pareschi@unife.it}).}\and Liu Liu\footnote{Department of Mathematics, The Chinese University of Hong Kong, Shatin, N.T., Hong Kong SAR ({\tt lliu@math.cuhk.edu.hk})}\and Lorenzo Pareschi$^*$\and Xueyu Zhu\footnote{Department of Mathematics, University of Iowa, Iowa City, IA 52242, USA. ({\tt xueyu-zhu@uiowa.edu.})}}

\maketitle
\begin{abstract}
Uncertainty in data is certainly one of the main problems in epidemiology, as shown by the recent COVID-19 pandemic. The need for efficient methods capable of quantifying uncertainty in the mathematical model is essential in order to produce realistic scenarios of the spread of infection. In this paper, we introduce a bi-fidelity approach to quantify uncertainty in spatially dependent epidemic models. The approach is based on evaluating a high-fidelity model on a small number of samples properly selected from a large number of evaluations of a low-fidelity model. In particular, we will consider the class of multiscale transport models recently introduced in^^>\cite{BDP, Bert3} as the high-fidelity reference and use simple two-velocity discrete models for low-fidelity evaluations. Both models share the same diffusive behavior and are solved with ad-hoc asymptotic-preserving numerical discretizations. A series of numerical experiments confirm the validity of the approach.
\end{abstract}


\tableofcontents

\section{Introduction}
Mathematical modeling in epidemiology has certainly experienced an impressive increase in recent times driven by the overwhelming effects of the COVID-19 pandemic. A large part of the research has been directed towards the construction of models capable of describing specific characteristics associated with the pandemic, in particular the presence of asymptomatic individuals, which were largely underestimated, especially in the early stages of the spread of the disease^^>\cite{Buonomo, Gatto, FBK2020,Peirlinck,Tang}. Another line of research was aimed at the construction of models capable of describing the spatial characteristics of the epidemic, in order to be able to properly assess the impact of containment measures, in particular with regard to the mobility on the territory and restrictions in areas with higher risk of infection^^>\cite{BDP, Veneziani2020, Veneziani2021, Colombo,SBKT2021,roques2020}. Other approaches took into account the network structure of connections^^>\cite{Bert}, social heterogeneity aspects^^>\cite{dimarco2021}, and the multiscale nature of the pandemic^^>\cite{bellomo2020multiscale}. See also^^>\cite{PuSa,LT21} and the survey^^>\cite{ABBDPTZ21} for some recent developments in epidemiological modelling based on the use of kinetic equations. 

Regardless of the characteristics of the model, a common aspect concerns the uncertainty of the data, which, especially in the first phase of the pandemic, largely underestimated the number of infected individuals. More precisely, the difficulty in correctly identifying all infected individuals in the early stages of the pandemic, due to structural limitations in performing large-scale screening and the inability to track the number of contacts, required the introduction of stochastic parameters into the models and the construction of related techniques for quantifying uncertainty^^>\cite{APZ, APZ2, Bert2, Bert3}.  

Among the various techniques of uncertainty quantification, the approaches based on stochastic strategies that do not necessarily require a-priori knowledge of the probability distribution of the uncertain parameters, as needed in the case of methods based on generalized polynomial chaos^^>\cite{Xiu}, are particularly interesting in view of a comparison with experimental data. On the other hand, the low convergence rate of Monte Carlo type sampling techniques poses serious limitations to their practical use.

In this context, multi-fidelity methods have shown to be able to efficiently alleviate such limitations through control variate techniques based on an appropriate use of low-fidelity surrogate models able to accelerate the convergence of stochastic sampling^^>\cite{ZNX14,NGX14,LPX21,DPM1,DPM2}. Specifically, general multi-fidelity approaches for kinetic equations have been developed in^^>\cite{DPM1, DPM2}, while bi-fidelity techniques with greedy sample selection in^^>\cite{LPX21,LX20}. We refer also to the recent survey in^^>\cite{DLPX21}. 

The present paper is devoted to extend the bi-fidelity method for transport equations in the diffusive limit developed in^^>\cite{LPX21} to the case of compartmental systems of multiscale equations designed to model mobility dynamics in an epidemic setting with uncertainty^^>\cite{BDP,Bert3}. For this purpose, the corresponding hyperbolic system recently introduced in^^>\cite{Bert,Bert2} will be used as a reduced low-fidelity model. The two models allow to correctly describe the hyperbolic dynamics of the movement of individuals over long distances together with the small-scale, high-density, diffusive nature typical of urban areas^^>\cite{BDP}. In addition, both models share the same diffusive limit in which it is possible to recover classical models of diffusive type^^>\cite{HS}.
From a numerical viewpoint, a space-time asymptotic-preserving discretization, which work uniformly in all regimes, has been adopted in combination with the bi-fidelity approach. This permits to obtain efficient stochastic asymptotic-preserving methods.    

The rest of the manuscript is organized as follows. In Section 2 we introduce the epidemic transport model with uncertainty in the simplified SIR compartmental setting, together with its diffusive limit. The corresponding reduced order two-velocity model used as low-fidelity surrogate is also discussed.
Next, we extend the previous modelling to more realistic compartmental settings based on the introduction of the exposed and asymptomatic compartments in Section 3. The details of the bi-fidelity method and the asymptotic-preserving IMEX Finite Volume scheme are then given in Section 4. Section 5 contains several numerical experiments that illustrate the performance of the bi-fidelity approach. Some conclusions are reported at the end of the manuscript.  

\section{Epidemic transport model with random inputs}
\label{section_SIR}
For simplicity, we first illustrate the modelling in the case of a classic SIR compartmental dynamic and subsequently we will extend our arguments to a more realistic SEIAR model, designed to take into account specific features of the COVID-19 pandemic, in Section \ref{section_SEIAR}.

\subsection{The high-fidelity epidemic transport model}
Let us consider a random vector ${\bf z}=(z_1,\ldots,z_{d})^T \in \RR^d$ characterizing possible sources of uncertainty due to the independent stochastic parameters $z_1,\ldots,z_{d}$, which may affect variables of the mathematical model, as well as parameters or initial conditions. Individuals at position $x\in\Omega\subset \RR$ at time $t$ moving with velocity $v \in [-1,1]$ are denote by $f_S=f_S(x,v,t,{\bf z})$, $f_I=f_I(x,v,t,{\bf z})$ and $f_R=f_R(x,v,t,{\bf z})$, which are the respective kinetic densities of susceptible $S$ (individuals who may be infected by the disease), infectious $I$ ( individuals who may transmit the disease) and removed $R$ (individuals healed or died due to the disease). We assume to have a population with subjects having no prior immunity and neglect the vital dynamics represented by births and deaths due to the time scale considered. 

The kinetic distribution is then given by
\[
f(x,v,t,{\bf z})=f_S(x,v,t,{\bf z})+f_I(x,v,t,{\bf z})+f_R(x,v,t,{\bf z}),
\]
and we recover their total density by integration over the velocity space 
\[
\rho(x,t,{\bf z})=\int_{-1}^{1} f(x,v_*,t,{\bf z})\,dv_*.
\]
As a consequence,
\begin{eqnarray}
S(x,t,{\bf z})&=&\int_{-1}^{1}  f_S(x,v,t,{\bf z})\,dv\nonumber\\
I(x,t,{\bf z})&=&\int_{-1}^{1}  f_I(x,v,t,{\bf z})\,dv \label{def.densitiesSIR}\\
R(x,t,{\bf z})&=&\int_{-1}^{1}  f_R(x,v,t,{\bf z})\,dv,\nonumber
\end{eqnarray}
with $\rho(x,t,{\bf z})=S(x,t,{\bf z})+I(x,t,{\bf z})+R(x,t,{\bf z})$, 
denote the density fractions of the population at position $x$ and time $t>0$ that are susceptible, infected and recovered respectively. In this setting, the kinetic densities satisfy the epidemic transport equations^^>\cite{BDP}
\begin{eqnarray}
\label{eqn:kinetic}
\nonumber
\frac{\partial f_S}{\partial t} + v_S \frac{\partial f_S}{\partial x} &=& -F(f_S, I) +\frac1{\tau_S}\left(\frac{S}{2}-f_S\right)\\
\label{eq:kineticc}
\frac{\partial f_I}{\partial t} + v_I \frac{\partial f_I}{\partial x} &=&  F(f_S, I)-\gamma f_I+\frac1{\tau_I}\left(\frac{I}{2}-f_I\right)\\
\nonumber
\frac{\partial f_R}{\partial t} + v_R \frac{\partial f_R}{\partial x} &=& \gamma f_I+\frac1{\tau_R}\left(\frac{R}{2}-f_R\right),
\end{eqnarray}
where 
$v_S=\vs(x) v$, $v_I=\vi(x) v$, $v_R=\vr(x) v$, $\vs,\vi,\vr \geq 0$ take into account the heterogeneities of geographical areas, and are thus chosen dependent on the spatial location. Similarly, also the relaxation times $\tau_S=\tau_S(x)$, $\tau_I=\tau_I(x)$ and $\tau_R=\tau_R(x)$. The quantity $\gamma=\gamma(x,{\bf z})$ is the recovery rate of infected, while the transmission of the infection is governed by an incidence function $F(\cdot,I)$ modeling the transmission of the disease^^>\cite{HWH00}
\begin{equation}
F(g,I)=\beta \frac{g I^p}{1+\kappa I},
\label{eq:incf}
\end{equation}
with the classic bi-linear case corresponding to $p = 1$, $\kappa=0$, even though it has been observed that an incidence rate that increases more than linearly with respect to the number of infected $I$ can occur under certain circumstances^^>\cite{CS78,BCV13,KM05}. The parameter $\beta=\beta(x,t,{\bf z})$ characterizes the average number of contacts per person per time, multiplied by the probability of disease transmission in a contact between a susceptible and an infectious subject; whereas $\kappa=\kappa(x,t,{\bf z}) > 0$ acts as an incidence damping parameter based on the self-protective behavior of the individual that arises from awareness of the epidemic risk as the disease progresses^^>\cite{Wang2020,Franco2020}. Notice that both $\beta$ and $\kappa$ may vary in time as a consequence of governmental control actions, such as mandatory wearing of masks or full lockdowns, and the increasing awareness of the epidemic risks among the population.

The standard threshold of epidemic models is the well-known reproduction number $R_0$, which defines the average number of secondary infections produced when one infected individual is introduced into a host population in which everyone is susceptible^^>\cite{HWH00}. This number determines when an infection can invade and persist in a new host population. For many deterministic infectious disease models, an infection begins in a fully susceptible population if and only if $R_0 > 1$. 
Assuming no inflow/outflow boundary conditions in $\Omega$, integrating over velocity/space and summing up the second equation in \eqref{eq:kineticc} we have
\[
\frac{\partial}{\partial t} \int_{\Omega} I(x,t,{\bf z})\,dx =  \int_{\Omega} F(S,I)\,dx-\int_{\Omega} \gamma(x,{\bf z}) I(x,t,{\bf z})\,dx \geq 0
\]
when
\begin{equation}
R_0(t,{\bf z})=\frac{\int_{\Omega} F(S,I)\,dx}{\int_{\Omega} \gamma(x,{\bf z}) I(x,t,{\bf z})\,dx} \geq 1.
\label{eq:R0}
\end{equation}
The above quantity, therefore, defines the stochastic reproduction number for system \eqref{eq:kineticc} describing the space averaged instantaneous variation of the number of infectious individuals at time $t>0$. This definition naturally extends locally by integrating over any subset of the computational domain $\Omega$ if one ignores the boundary flows.


\subsubsection{Macroscopic formulation and diffusion limit}
\label{sect:diff-1Dkinetic}
Let us introduce the flux functions
\begin{equation}
J_S={\lambda_S} \int_{-1}^{1}  v f_S(x,v,t,{\bf z})\,dv,\quad J_I={\lambda_I}\int_{-1}^{1}  v f_I(x,v,t,{\bf z})\,dv,\quad J_R={\lambda_R}\int_{-1}^{1}  v f_R(x,v,t,{\bf z})\,dv.
\label{eq.fluxes}
\end{equation}
Then, integrating the system \eqref{eq:kineticc} against $v$, it is straightforward to get the following set of equations for the macroscopic densities of commuters
\begin{eqnarray}
\nonumber
\frac{\partial S}{\partial t} + \frac{\partial J_S}{\partial x} &=& -F(S, I)\\
\label{eq:density}
\frac{\partial I}{\partial t} + \frac{\partial J_I}{\partial x} &=& F(S, I) -\gamma I\\
\nonumber
\frac{\partial R}{\partial t} + \frac{\partial J_R}{\partial x} &=& \gamma I
\end{eqnarray}
whereas the flux functions satisfy
\begin{eqnarray}
\nonumber
\frac{\partial J_S}{\partial t} +  {\vs^2} \int_{-1}^1 v^2 \frac{\partial f_S}{\partial x}\,dv &=& -F(J_S, I)-\frac{J_S}{\tau_S} \\
\label{eq:flux}
\frac{\partial J_I}{\partial t} +  {\vi^2} \int_{-1}^1 v^2 \frac{\partial f_I}{\partial x}\,dv &=& -\frac{\lambda_I}{\lambda_S}F(J_S, I) - \gamma J_I-\frac{J_I}{\tau_I} \\
\nonumber
\frac{\partial J_R}{\partial t} +  {\vr^2} \int_{-1}^1 v^2 \frac{\partial f_R}{\partial x}\,dv &=& -\frac{\lambda_R}{\lambda_I} \gamma J_I-\frac{J_R}{\tau_R} .
\end{eqnarray}
By introducing the space dependent diffusion coefficients 
\begin{equation}
D_S=\frac13\lambda_S^2\tau_S,\quad D_I=\frac13\lambda_I^2\tau_I,\quad D_R=\frac13\lambda_R^2\tau_R,
\label{eq:diffcf}
\end{equation}
which characterize the diffusive transport mechanism of susceptible, infectious and removed, respectively,
and keeping the above quantities fixed while letting the relaxation times $\tau_{S,I,R}$ to zero, we get from the r.h.s. in \eqref{eq:kineticc}
\[
\begin{split}
&f_S=S/2,\quad f_I=I/2 ,\quad f_R=R/2,
\end{split}
\]
and consequently, from \eqref{eq:flux}, we obtain a proportionality relation between the fluxes and the spatial derivatives (Fick's law):
\begin{equation}
J_S = -D_S \frac{\partial S}{\partial x},\quad  J_I = -D_I\frac{\partial I}{\partial x},\quad J_R = -D_R\frac{\partial R}{\partial x}.
\label{eq:fick}
\end{equation}
Thus, substituting \eqref{eq:fick} into \eqref{eq:density} we get the diffusion system^^>\cite{MWW, Sun, Webb,HS}
\begin{eqnarray}
\nonumber
\frac{\partial S}{\partial t} &=&  -F(S, I)+\frac{\partial}{\partial x} \left({D_S}\frac{\partial S}{\partial x}\right)\\
\label{eq:diff}
\frac{\partial I}{\partial t} &=&  F(S, I)-\gamma I+\frac{\partial}{\partial x} \left({D_I}\frac{\partial I}{\partial x}\right)\\
\nonumber
\frac{\partial R}{\partial t} &=&  \gamma I+\frac{\partial}{\partial x}\left({D_R}\frac{\partial R}{\partial x}\right).
\end{eqnarray}
We remark that the model's capability to account for different regimes, hyperbolic or parabolic, accordingly to the space dependent values $\tau_S$, $\tau_I$, $\tau_R$, makes it suitable for describing the dynamics of populations composed of human beings. Indeed, it is clear that the daily routine is a complex mixing of individuals moving at the scale of a city and individuals moving among different urban centers. 
In this situation, it seems reasonable to avoid, due to the lack of microscopic information and the high complexity, the description of the details of movements within an urban area and to describe this aspect through a diffusion operator. On the other hand, individuals when moving from one city to another follow well established connections for which a hyperbolic setting is certainly more  appropriate.

\subsection{The low-fidelity two-velocity epidemic model}
The low-fidelity model, is based on considering individuals moving in two opposite directions (indicated by signs ``+'' and ``-''), with velocities $\pm \lambda_S$ for susceptible, $\pm \lambda_I$ for infectious and $\pm \lambda_R$ for removed, we can describe the spatio-temporal dynamics of the population through the following two-velocity epidemic model^^>\cite{Bert}
\begin{eqnarray}
	\frac{\partial S^{\pm}}{\partial t} + \lambda_S \frac{\partial S^{\pm}}{\partial x} &=& - F(S^{\pm},I) \mp \frac{1}{2\tau_S}\left(S^+ - S^-\right)				
	\nonumber\\ 
	\frac{\partial I^{\pm}}{\partial t} + \lambda_I \frac{\partial I^{\pm}}{\partial x} &=&  F(S^{\pm},I) -\gamma I^{\pm} \mp \frac{1}{2\tau_I}\left(I^+ - I^-\right)								
	\label{eq.SIR_kinetic_diag}	\\ 
	\frac{\partial R^{\pm}}{\partial t} + \lambda_R \frac{\partial R^{\pm}}{\partial x} &=& \gamma I^{\pm} \mp \frac{1}{2\tau_R}\left(R^+ - R^-\right).									
    \nonumber
\end{eqnarray}
In the above system, individuals $S(x,t,{\bf z})$, $I(x,t,{\bf z})$ and $R(x,t,{\bf z})$ are defined as
\begin{equation*}
 	S = S^+ + S^- , \quad 
 	I = I^+ + I^- , \quad 
 	R = R^+ + R^- .	
\end{equation*}
The transmission of the infection is governed by the same incidence function as in the high-fidelity model, defined in \eqref{eq:incf}. Also the definition of the basic reproduction number $R_0$ results the same previously introduced in \eqref{eq:R0}.

\subsubsection{Macroscopic formulation and diffusion limit}
\label{sect:diff-LF}
If we now introduce the fluxes, defined by
\begin{equation}
 	J_S = \lambda_S \left(S^+ - S^-\right) , \quad 
 	J_I = \lambda_I \left(I^+ - I^-\right) , \quad 
 	J_R = \lambda_R \left(R^+ - R^-\right) ,	  	\label{eq.fluxes_kinetic_diag}
\end{equation}
we obtain a hyperbolic model equivalent to \eqref{eq.SIR_kinetic_diag}, but presenting a macroscopic description of the propagation of the epidemic at finite speeds, for which the densities follow system \eqref{eq:density}
and equations of fluxes read 
\begin{eqnarray}
	\frac{\partial J_S}{\partial t} + \lambda_S^2 \frac{\partial S}{\partial x} &=& -F(J_S,I)  -\frac{J_S}{\tau_S} 
	\nonumber\\ 
	\frac{\partial J_I}{\partial t} + \lambda_I^2 \frac{\partial I}{\partial x} &=& \frac{\lambda_I}{\lambda_S}F(J_S,I)  -\gamma J_I -\frac{J_I}{\tau_I}	
	\label{eq.SIR_kinetic-fluxes}\\ 
	\frac{\partial J_R}{\partial t} + \lambda_R^2 \frac{\partial R}{\partial x} &=& \frac{\lambda_R}{\lambda_I}\gamma J_I -\frac{J_R}{\tau_R} .
	\nonumber
\end{eqnarray}

Let us now consider the behavior of the low-fidelity model in diffusive regimes. To this aim, we introduce the diffusion coefficients 
\begin{equation}
D_S=\lambda_S^2 \tau_S,\quad D_I=\lambda_I^2 \tau_I,\quad D_R=\lambda_R^2 \tau_R.
\label{eq:diff2}
\end{equation}
As for the previous model, the diffusion limit of the system is formally recovered letting the relaxation times $\tau_{S,I,R}\to 0$, 
while keeping the diffusion coefficients \eqref{eq:diff2} finite. Under this scaling, from \eqref{eq.SIR_kinetic-fluxes} we get equations \eqref{eq:fick},
which inserted into \eqref{eq:density} lead again to the parabolic reaction-diffusion system \eqref{eq:diff}.

\begin{remark}
\label{remark_consistency}
An important aspect in the bi-fidelity approach here proposed, is that the high-fidelity model and the low-fidelity one exactly coincide in the diffusive limit. The only difference lays in the definition of the diffusion coefficients (see eqs. \eqref{eq:diffcf} and eqs.  \eqref{eq:diff2}). As a consequence, in such a regime the bi-fidelity method achieves the maximum accuracy.
\end{remark}

\section{Extension to more epidemic compartments}
\label{section_SEIAR}
To account for more realistic models to analyze the evolution of the ongoing COVID-19 pandemic, we consider extending the simple SIR compartmentalization by taking into account two additional population compartments, $E$ and $A$, resulting in a SEIAR model^^>\cite{Bert2,Bert3}. Subjects in the $E$ compartment are the exposed, hence infected but not yet infectious, being in the latent period. Moreover, among the infectious subjects, we distinguish the population between a group of individuals $I$ who will develop severe symptoms and a group of individuals $A$ who will never develop symptoms or, if they do, these will be very mild. This feature turns out to be essential to correctly analyze the evolution of COVID-19. In fact, it has been shown that individuals belonging to the $A$ group are very difficult to detect and isolate, contributing more strongly to the spread of the virus than the more easily detectable $I$ individuals^^>\cite{Gatto,Peirlinck,Tang}.

\subsection{The high-fidelity SEIAR transport model}
Let us now consider a population at position $x\in\Omega$ moving with velocity directions $v \in [-1,1]$, still taking into account possible uncertainties related to the random vector ${\bf z}$. Defining the kinetic densities of susceptible $f_S=f_S(x,v,t,{\bf z})$, exposed $f_E=f_E(x,v,t,{\bf z})$, severe symptomatic infected $f_I=f_I(x,v,t,{\bf z})$, mildly symptomatic or asymptomatic infected $f_A=f_A(x,v,t,{\bf z})$ and removed (healed or deceased) $f_R=f_R(x,v,t,{\bf z})$, the kinetic distribution of the population results
\[
f(x,v,t,{\bf z})=f_S(x,v,t,{\bf z})+f_E(x,v,t,{\bf z})+f_I(x,v,t,{\bf z})+f_A(x,v,t,{\bf z})+f_R(x,v,t,{\bf z}).
\]
In addition to \eqref{def.densitiesSIR}, we have that 
\begin{equation*}
E(x,t,{\bf z})=\int_{-1}^{1}  f_E(x,v,t,{\bf z})\,dv,\qquad A(x,t,{\bf z})=\int_{-1}^{1}  f_A(x,v,t,{\bf z})\,dv.
\end{equation*}
In this setting, the kinetic densities satisfy the following transport equations^^>\cite{Bert3}
\begin{eqnarray}
\nonumber
\frac{\partial f_S}{\partial t} + v_S \frac{\partial f_S}{\partial x} &=& -F(f_S, I)-F_A(f_S, A) +\frac1{\tau_S}\left(\frac{S}{2}-f_S\right)\\
\nonumber
\frac{\partial f_E}{\partial t} + v_E \frac{\partial f_E}{\partial x} &=&  F(f_S, I)+F_A(f_S, A)-a f_E+\frac1{\tau_E}\left(\frac{E}{2}-f_E\right)\\
\label{eq:kineticc-SEIAR}
\frac{\partial f_I}{\partial t} + v_I \frac{\partial f_I}{\partial x} &=& a\sigma f_E -\gamma_I f_I+\frac1{\tau_I}\left(\frac{I}{2}-f_I\right)\\
\nonumber
\frac{\partial f_A}{\partial t} + v_A \frac{\partial f_A}{\partial x} &=& a(1-\sigma) f_E -\gamma_A f_A+\frac1{\tau_A}\left(\frac{A}{2}-f_A\right)\\
\nonumber
\frac{\partial f_R}{\partial t} + v_R \frac{\partial f_R}{\partial x} &=& \gamma_I f_I+\gamma_A f_A+\frac1{\tau_R}\left(\frac{R}{2}-f_R\right)
\end{eqnarray}
The quantities $\gamma_I=\gamma_I(x,{\bf z})$ and $\gamma_A=\gamma_A(x,{\bf z})$ are the recovery rates of symptomatic and asymptomatic infected (inverse of the infectious periods), respectively, while $a=a(x,{\bf z})$ represents the inverse of the latency period and $\sigma=\sigma(x,{\bf z})$ is the probability rate of developing severe symptoms^^>\cite{Tang,Gatto,Buonomo}.
In this model, the transmission of the infection is governed by two different incidence functions, $F(\cdot,I)$ and $F_A(\cdot,A)$, simply to distinguish between the behavior of $I$ and $A$ individuals. Analogously to \eqref{eq:incf},
\begin{equation}
F_A(g,A)=\beta_A \frac{g A^p}{1+\kappa_A A},
\label{eq:incf-SEIAR}
\end{equation}
where a different contact rate, $\beta_A = \beta_A(x,t,{\bf z})$, and coefficient $\kappa_A=\kappa_A(x,t,{\bf z})$ are taken into account for mildly/no symptomatic people. For the derivation of the reproduction number $R_0$ of this SEIAR kinetic model, which results
\begin{equation}
\begin{aligned}
R_0 (t,{\bf z}) &= \frac{\int_\Omega F_I(S,I) \,dx}{\int_\Omega \gamma_I(x,{\bf z}) I(x,t,{\bf z}) \,dx} \cdot \frac{\int_\Omega a(x,{\bf z})\sigma (x,{\bf z}) E(x,t,{\bf z}) \,dx}{\int_\Omega a(x,{\bf z}) E(x,t,{\bf z}) \,dx} \\&+ \frac{\int_\Omega F_A(S,A) \,dx}{\int_\Omega \gamma_A(x,{\bf z}) A(x,t,{\bf z})\, dx} \cdot \frac{\int_\Omega a(x,{\bf z})(1-\sigma(x,{\bf z})) E(x,t,{\bf z}) \,dx}{\int_\Omega a(x,{\bf z}) E(x,t,{\bf z})\, dx} \,,
\end{aligned}
\label{eq.R0_2}
\end{equation}
the reader is invited to refer to^^>\cite{Bert2,Bert3}.

When introducing the same definition \eqref{eq.fluxes} of flux for the additional compartments, $J_E$ and $J_A$, integrating system \eqref{eq:kineticc} in $v$, we get the following set of equations for the macroscopic densities
\begin{eqnarray}
\nonumber
\frac{\partial S}{\partial t} + \frac{\partial J_S}{\partial x} &=& -F(S, I)-F_A(S, A)\\
\nonumber
\frac{\partial E}{\partial t} + \frac{\partial J_E}{\partial x} &=& F(S, I)+F_A(S, A)-aE\\
\label{eq:density-1Dkinetic-SEIAR}
\frac{\partial I}{\partial t} + \frac{\partial J_I}{\partial x} &=& a\sigma E -\gamma_I I\\
\nonumber
\frac{\partial A}{\partial t} + \frac{\partial J_A}{\partial x} &=& a(1-\sigma) E -\gamma_A A\\
\nonumber
\frac{\partial R}{\partial t} + \frac{\partial J_R}{\partial x} &=& \gamma_I I+\gamma_A A
\end{eqnarray}
and for the fluxes
\begin{eqnarray}
\nonumber
\frac{\partial J_S}{\partial t} +  {\vs^2} \int_{-1}^1 v^2 \frac{\partial f_S}{\partial x}\,dv &=& -F(J_S, I)-F_A(J_S,A)-\frac{J_S}{\tau_S} \\
\nonumber
\frac{\partial J_E}{\partial t} +  {\ve^2} \int_{-1}^1 v^2 \frac{\partial f_E}{\partial x}\,dv &=& \frac{\lambda_E}{\lambda_S}\left(F(J_S, I)+F_A(J_S,A)\right)-a J_E -\frac{J_E}{\tau_E} \\
\label{eq:flux-1Dkinetic-SEIAR}
\frac{\partial J_I}{\partial t} +  {\vi^2} \int_{-1}^1 v^2 \frac{\partial f_I}{\partial x}\,dv &=& \frac{\lambda_I}{\lambda_E}a\sigma J_E - \gamma_I J_I-\frac{J_I}{\tau_I} \\
\nonumber
\frac{\partial J_A}{\partial t} +  {\va^2} \int_{-1}^1 v^2 \frac{\partial f_A}{\partial x}\,dv &=& \frac{\lambda_A}{\lambda_E}a(1-\sigma) J_E - \gamma_A J_A-\frac{J_A}{\tau_A} \\
\nonumber
\frac{\partial J_R}{\partial t} + {\vr^2} \int_{-1}^1 v^2 \frac{\partial f_R}{\partial x}\,dv &=& \frac{\lambda_R}{\lambda_I} \gamma_I J_I+\frac{\lambda_R}{\lambda_A} \gamma_A J_A
-\frac{J_R}{\tau_R} .
\end{eqnarray}

Moreover, defining also $D_E=\frac13\lambda_E^2\tau_E$ and $D_A=\frac13\lambda_A^2\tau_A$ and following the same procedure discussed in Section \ref{sect:diff-1Dkinetic}, we recover this SEIAR system in the diffusive regime:
\begin{eqnarray}
\nonumber
\frac{\partial S}{\partial t} &=& -F(S, I) -F_A(S, A) + \frac{\partial}{\partial x} \left({D_S}\frac{\partial S}{\partial x}\right) \\
\nonumber
\frac{\partial E}{\partial t} &=& F(S, I) + F_A(S, A) -a E +\frac{\partial}{\partial x} \left({D_E}\frac{\partial E}{\partial x}\right) \\
\label{eq:diff-SEIAR}
\frac{\partial I}{\partial t}  &=&  a\sigma E-\gamma_I I +\frac{\partial}{\partial x} \left({D_I}\frac{\partial I}{\partial x}\right)\\
\nonumber
\frac{\partial A}{\partial t}  &=&  a(1-\sigma) E-\gamma_A A +\frac{\partial}{\partial x} \left({D_A}\frac{\partial A}{\partial x}\right)\\
\nonumber
\frac{\partial R}{\partial t}  &=& \gamma_I I +\gamma_A A +\frac{\partial}{\partial x} \left({D_R}\frac{\partial R}{\partial x}\right) .
\end{eqnarray}


\subsection{The low-fidelity SEIAR transport model}
The low-fidelity model is obtained, also with the more complex SEIAR compartmentalization, considering the discrete-velocity case with only two opposite velocities $\pm \lambda$^^>\cite{Bert3}:
\begin{eqnarray}
\frac{\partial S^{\pm}}{\partial t} \pm \lambda_S \frac{\partial S^{\pm}}{\partial x} &=& -F(S^{\pm}, I) -F_A(S^{\pm}, A) + \frac{1}{2\tau_S}\left(S^\mp - S^\pm\right)
\nonumber\\
\frac{\partial E^{\pm}}{\partial t} \pm \lambda_E \frac{\partial E^{\pm}}{\partial x} &=& F(S^{\pm}, I) +F_A(S^{\pm},A) -a E^{\pm} + \frac{1}{2\tau_E}\left(E^\mp - E^\pm\right)
\nonumber\\
\frac{\partial I^{\pm}}{\partial t} \pm \lambda_I \frac{\partial I^{\pm}}{\partial x} &=& a \sigma E^{\pm} -\gamma_I I^{\pm} + \frac1{2\tau_I}\left(I^\mp - I^\pm\right)
\label{eq.SEIARkinetic}\\
\frac{\partial A^{\pm}}{\partial t} \pm \lambda_A \frac{\partial A^{\pm}}{\partial x} &=& a(1-\sigma) E^{\pm} -\gamma_A A^{\pm} + \frac{1}{2\tau_A}\left(A^\mp - A^\pm\right)
\nonumber\\
\frac{\partial R^{\pm}}{\partial t} \pm \lambda_R \frac{\partial R^{\pm}}{\partial x} &=& \gamma_I I^{\pm} + \gamma_A A^{\pm} + \frac1{2\tau_R}\left(R^\mp - R^\pm\right)\, .
\nonumber
\end{eqnarray}

When defining fluxes $J_E$ and $J_A$ as in \eqref{eq.fluxes_kinetic_diag}, the equivalent hyperbolic model underlying the macroscopic formulation of the spatial propagation of an epidemic is obtained. The evolution of the densities follows \eqref{eq:density-1Dkinetic-SEIAR},
while for the fluxes we have
\begin{eqnarray}
\frac{\partial J_S}{\partial t} + \lambda_S^2 \frac{\partial S}{\partial x} &=& -F(J_S, I) -F_A(J_S, A) - \frac{J_S}{\tau_S} 
\nonumber\\
\frac{\partial J_E}{\partial t} + \lambda_E^2 \frac{\partial E}{\partial x} &=& \frac{\lambda_E}{\lambda_S}\left(F(J_S, I) + F_A(J_S, A)\right) - a J_E - \frac{J_E}{\tau_E} 
\nonumber\\
\frac{\partial J_I}{\partial t} + \lambda_I^2 \frac{\partial I}{\partial x} &=&  \frac{\lambda_I}{\lambda_E} a \sigma J_E  - \gamma_I J_I-\frac{J_I}{\tau_I} 
\label{eq.SEIARmacro-fluxes}\\
\frac{\partial J_A}{\partial t} + \lambda_A^2 \frac{\partial A}{\partial x} &=&  \frac{\lambda_A}{\lambda_E} a (1- \sigma) J_E  - \gamma_A J_A-\frac{J_A}{\tau_A} \nonumber\\
\frac{\partial J_R}{\partial t} + \lambda_R^2 \frac{\partial R}{\partial x} &=& \frac{\lambda_R}{\lambda_I} \gamma_I J_I + \frac{\lambda_R}{\lambda_A} \gamma_A J_A -\frac{J_R}{\tau_R}  .
\nonumber
\end{eqnarray}

As for the previous cases, the diffusion limit of the system is formally recovered letting the relaxation times $\tau_{S,E,I,A,R}\to 0$, while keeping the diffusion coefficients, $D_i=\lambda_i^2 \tau_i, i\in\{S,E,I,A,R\}$, finite. Thus, the same procedure presented for the SIR compartmentalization in Section \ref{sect:diff-LF} lead to the SEIAR parabolic system \eqref{eq:diff-SEIAR}.

\section{An asymptotic-preserving bi-fidelity numerical method}
In this Section, we present the details of the numerical method used to solve the stochastic problem following the bi-fidelity asymptotic-preserving scheme proposed in^^>\cite{NGX14, ZNX14}. 
For the high-fidelity model, the numerical scheme is structured in agreement with a discrete ordinate method in velocity with the even and odd parity formulation^^>\cite{DP,JPT}. Both the high-fidelity and low-fidelity solvers use Finite Volume Method (FVM) in space and achieve asymptotic preservation in time using suitable IMEX Runge-Kutta schemes^^>\cite{Bos1, Bos2}. This permits to obtain a numerical scheme able to deal with the diffusion limit of the mathematical models without loosing consistency, and for which the time step size of the temporal discretization is not subject to excessive restrictions related to the smallness of the scaling parameters $\tau_i$. 
For simplicity, we illustrate the numerical method in the case of the simpler SIR model. The application of the same numerical scheme results straightforward for the case of the SEIAR compartmentalization.

\subsection{Asymptotic-preserving IMEX Finite Volume scheme}
\label{sect_AP_IMEX}
In this Section, we present the AP-IMEX Finite Volume scheme adopted to solve the SIR model at each stochastic collocation point selected for the bi-fidelity approximation.

The asymptotic-preserving IMEX method and the corresponding even and odd parities formulation was introduced in^^>\cite{BDP} for an SIR kinetic transport model in 2D domains. 
According to^^>\cite{JPT,GJL}, for $v>0$, we can define the even and odd parities for the high-fidelity SIR kinetic transport model \eqref{eqn:kinetic} as follows: 
\begin{eqnarray}
\nonumber
r_S(v) = \frac{1}{2}\left(f_S(v) + f_S(-v)\right), & \qquad 
j_S(v) = \frac{\lambda_S}{2}\left(f_S(v) - f_S(-v)\right), \\[4pt]
\nonumber
r_I(v) = \frac{1}{2}\left(f_I(v) + f_I(-v)\right), & \qquad 
j_I(v) = \frac{\lambda_I}{2}\left(f_I(v) - f_I(-v)\right), \\[4pt]
\nonumber
r_R(v) = \frac{1}{2}\left(f_R(v) + f_R(-v)\right), & \qquad 
j_R(v) = \frac{\lambda_R}{2}\left(f_R(v) - f_R(-v)\right). 
\end{eqnarray}
An equivalent formulation of \eqref{eqn:kinetic} can be written as
\begin{equation}
\label{eq:AP_Kinetic}
\begin{split}
\frac{\partial r_S}{\partial t} + v\frac{\partial j_S}{\partial x} & =  - F(r_S, I) + \frac{1}{\tau_S}\left(\frac{1}{2}S - r_S\right) \\
\frac{\partial r_I}{\partial t} + v\frac{\partial j_I}{\partial x} & =  F(r_S, I) - \gamma r_I + \frac{1}{\tau_I}\left(\frac{1}{2}I - r_I\right) \\
\frac{\partial r_R}{\partial t} + v\frac{\partial j_R}{\partial x} & =  \gamma r_I + \frac{1}{\tau_R}\left(\frac{1}{2}R - r_R\right) \\
\frac{\partial j_S}{\partial t} + \lambda_S^2\, v\frac{\partial r_S}{\partial x} & =  - F(j_S, I)- \frac{1}{\tau_S}j_S \\
\frac{\partial j_I}{\partial t} + \lambda_I^2\, v\frac{\partial r_I}{\partial x} & =  F(j_S, I) - \gamma j_I - \frac{1}{\tau_I}j_I \\
\frac{\partial j_R}{\partial t} + \lambda_R^2\, v\frac{\partial r_R}{\partial x} & =  \gamma j_I - \frac{1}{\tau_R}j_R\, ,
\end{split}
\end{equation}
where
\begin{equation}
S=2\int_{0}^1  r_S\,dv,\quad
I=2\int_{0}^1  r_I\,dv,\quad
R=2\int_{0}^1  r_R\,dv.
\label{eq:dv}
\end{equation}
The above densities can be approximated by a Gauss-Legendre quadrature rule. This leads to a discrete velocity setting, usually referred to as the discrete ordinate method, where we approximate
\begin{equation*}
S \approx S_M = \sum_{i=1}^{N_G} w_i \,r_S(\zeta_i)\,\quad 
I \approx I_M = \sum_{i=1}^{N_G} w_i \,r_I(\zeta_i)\, \quad
R \approx R_M = \sum_{i=1}^{N_G} w_i \,r_R(\zeta_i)\,
\end{equation*}
where $w_i$ and $\zeta_i$ are the $N_G$ standard Gauss-Legendre quadrature weights and points in $[-1,1]$, and $N_G=N_v$, number of chosen discrete velocities.

Assuming for simplicity of notation that $\tau_{S,I,R} = \tau$, we can write the above system in the following compact form
\begin{equation}
\label{eq:AP_Kinetic_compact}
\begin{split}
\frac{\partial \boldsymbol{r}}{\partial t} + v \frac{\partial \boldsymbol{j}}{\partial x} &= \boldsymbol{E}(\boldsymbol{r}) - \frac{1}{\tau}\left( \boldsymbol{r} - \frac{\boldsymbol{R}}{2}\right)\\
\frac{\partial \boldsymbol{j}}{\partial t} + \boldsymbol{\Lambda}^2 v \frac{\partial \boldsymbol{r}}{\partial x} &= \boldsymbol{E}(\boldsymbol{j}) - \frac{1}{\tau}\boldsymbol{j}\,,
\end{split}
\end{equation}
with
\[
\boldsymbol{r} = \left( r_S,r_I,r_R \right)^T, \quad
\boldsymbol{j} = \left( j_S,j_I,j_R \right)^T, \quad
\boldsymbol{R} = \left( S,I,R \right)^T,\quad
\boldsymbol{\Lambda} = \mathrm{diag}\{\lambda_S,\lambda_I,\lambda_R\},
\]
\[
\boldsymbol{E}(\boldsymbol{r}) = \left( - F(r_S, I), F(r_S, I)- \gamma r_I, \gamma r_I \right)^T ,\quad
\boldsymbol{E}(\boldsymbol{j})\left( - F(j_S, I), F(j_S, I)- \gamma j_I, \gamma j_I  \right)^T.
\]
Following^^>\cite{Bos2}, the Implicit-Explicit (IMEX) Runge-Kutta discretization that we consider for system \eqref{eq:AP_Kinetic_compact} consists in computing the internal stages
\begin{equation}
\begin{aligned}
\boldsymbol{r} ^{(k)} &= \boldsymbol{r} ^n -  \Delta t \sum_{j=1}^{k} a_{kj} \left(v\frac{\partial \boldsymbol{j} ^{(j)}}{\partial x} +\frac1{\tau}\left(\boldsymbol{r}^{(j)} - \frac{\boldsymbol{R}^{(j)}}{2}\right)\right) + \Delta t \sum_{j=1}^{k-1} \tilde{a}_{kj} \boldsymbol{E}\left(\boldsymbol{r}^{(j)}\right)
\\
\boldsymbol{j}^{(k)} &= \boldsymbol{j}^n -  \Delta t \sum_{j=1}^{k-1} \tilde{a}_{kj} \left(\boldsymbol{\Lambda}^2 v \frac{\partial \boldsymbol{r}^{(j)}}{\partial x} -\boldsymbol{E}\left(\boldsymbol{j}^{(j)}\right)\right) 
- \Delta t \sum_{j=1}^{k} a_{kj} \frac1{\tau}\boldsymbol{j}^{(j)},
\end{aligned}
\label{eq.iterIMEX}
\end{equation}
and evaluating the final numerical solution
\begin{equation}
\begin{aligned}
\boldsymbol{r} ^{n+1} &= \boldsymbol{r} ^n -  \Delta t \sum_{k=1}^{s} b_{k} \left(v\frac{\partial \boldsymbol{j} ^{(k)}}{\partial x} +\frac1{\tau}\left(\boldsymbol{r}^{(k)} - \frac{\boldsymbol{R}^{(k)}}{2}\right)\right) + \Delta t \sum_{k=1}^{s} \tilde{b}_{k} \boldsymbol{E}\left(\boldsymbol{r}^{(k)}\right)
\\
\boldsymbol{j}^{n+1} &= \boldsymbol{j}^n -  \Delta t \sum_{k=1}^{s} \tilde{b}_{k} \left(\boldsymbol{\Lambda}^2 v \frac{\partial \boldsymbol{r}^{(k)}}{\partial x} -\boldsymbol{E}\left(\boldsymbol{j}^{(k)}\right)\right) 
- \Delta t \sum_{k=1}^{s} b_{k} \frac1{\tau}\boldsymbol{j}^{(k)}.
\end{aligned}
\label{eq.finalIMEX}
\end{equation}
To properly compute the implicit terms $\frac{\boldsymbol{R}^{(j)}}{2}$ at each Runge-Kutta internal step explicitly, we refer to^^>\cite{JHL}.
Matrices $\tilde A = (\tilde a_{kj})$, with $\tilde a_{kj} = 0 $ for $ j\geq k$, and $A = (a_{kj})$, with $a_{kj} = 0 $ for $ j > k$ are $s \times s$ matrices, with $s$ number of Runge-Kutta stages, defining respectively the explicit and the implicit part of the scheme, and vectors $\tilde b = (\tilde b_1, ...,\tilde b_s)^T$ and $b = (b_1, ...,b_s)^T$ are the quadrature weights.
Furthermore, referring to^^>\cite{Bos1, Bos2}, if the following relations hold, 
\begin{equation*}
a_{kj} = b_j, \qquad j = 1,\ldots,s ,\qquad
\tilde a_{kj} = \tilde b_j, \qquad j = 1,\ldots,s-1, 
\end{equation*}
the method is said to be globally stiffly accurate (GSA).
It is worth to notice that this definition states also that the numerical solution of a GSA IMEX Runge-Kutta scheme coincides exactly with the last internal stage of the scheme. 
Since the GSA property is fundamental to preserve the correct diffusion limit and to achieve asymptotic-preservation stability in stiff regimes^^>\cite{BDP,Bert}, in the sequel, the GSA BPR(4,4,2) scheme presented in^^>\cite{Bos2} is chosen, characterized by $s=4$ stages for the implicit part, 4 stages for the explicit part and 2nd order of accuracy.

At each internal stage of the IMEX scheme \eqref{eq.iterIMEX}, we apply a Total-Variation-Diminishing (TVD) Finite Volume discretization to evaluate the numerical fluxes^^>\cite{Bert,Bert4}. 
To achieve second order accuracy also in space, while avoiding the occurrence of spurious oscillations, a classical minmod slope limiter has been adopted.

The same AP-IMEX Finite Volume scheme is adopted also to solve the low-fidelity SIR model \eqref{eq:density}-\eqref{eq.SIR_kinetic-fluxes}, as fully presented in^^>\cite{Bert}. 
The reader is invited to refer to^^>\cite{Bert,Bert2,Bert3,BDP} for further details on the properties of the chosen numerical scheme applied to epidemic models.

\subsection{Bi-fidelity stochastic collocation}
\label{Bifi approx}

In this Section, we briefly review the bi-fidelity method developed in^^>\cite{NGX14, ZNX14}. The methods make use of low-fidelity models to effectively inform the selection of  representative points in the parameter space and then employ this information to construct  accurate approximations to high-fidelity solutions.  
To facilitate future discussion,
we denote the expensive high-fidelity solution $u^H({\bf{z}})$  and the cheap low-fidelity solution $u^L({\bf{z}})$  
for any given random parameter ${\bf{z}}\in I_z \subset \RR^d$, where $I_z$ is the domain of the random parameter ${\bf{z}}$. 

The basic idea of the bi-fidelity approximation 
is to construct an inexpensive surrogate $u^B({\bf z})$ of  the high-fidelity solution 
in the following non-intrusive  manner:
\begin{equation}\label{UB_1}
 u^B({\bf{z}}) = \sum_{k=1}^nc_k({\bf{z}}) u^H({\bf{z}}_{k}),
\end{equation}
where $n$ is the  number of the selected parameter points in the parameter space. If $n$ is small and the coefficient $c_k({\bf{z}})$ can be efficiently and accurately approximated, an efficient bi-fidelity approximation can be constructed. Once the surrogate is constructed, the statistics of high-fidelity solutions can be quickly approximated by evaluating the bi-fidelity surrogate via Monte Carlo methods or quadrature rules.

There are two major questions key to the performance of the above bi-fidelity algorithms: (a) how to select the representative points ${{\bf{z}}}_k$ effectively? (b) how to  efficiently construct the bi-fidelity approximation for  any new given ${\bf{z}}$ but avoiding requiring a high-fidelity simulation? 
\bigskip

{\bf Subset Selection}.
Existing  predefined or structured nodes (e.g., sparse grids,  cubature rules, etc.) often
grows fast in high dimensions. Therefore,
these options cannot easily accommodate the current situation where we would like the size $n$ of $\gamma_n$ to be 
small and also arbitrary.
Alternatively, adaptive approaches based on high-fidelity models to explore the parameter space suit our needs. However, this typically requires a large number of high-fidelity samples, which might not be computationally affordable. In contrast, the low-fidelity model is inexpensive to evaluate and it mimics the variations of high-fidelity solutions in the parameter space. This motivates us to employ the inexpensive low-fidelity model to learn and explore
the behaviors of the high-fidelity model in the parameter space. 

We shall identify important points iteratively by a greedy approach^^>\cite{NGX14, ZNX14}.  
Specifically, we denote a candidate sample set $\Gamma_N=\{{\bf{z}}_1, {\bf{z}}_2,\hdots, {\bf{z}}_N\}$, which is assumed to be large enough to cover the parameter space $I_z$.  Initially, denote $\gamma_0=\{\}$ and assume we have the first $k$ important points  $\gamma_k=\{{\bf{z}}_{i_1}, {\bf{z}}_{i_2}, \hdots, {\bf{z}}_{i_k}\}$ available at the $k$-th iteration. Denote the snapshot matrix $u^L(\gamma_k)= \{u({\bf{z}})| {\bf{z}} \in \gamma_k\}$ and the corresponding spanned approximation space $U^L(\gamma_k) = \text{span}\{u^L({\bf{z}})|{\bf{z}}\in\gamma_k\}$. The corresponding high-fidelity approximation space can be defined similarly , $U^H(\gamma_n) = \text{span}\{u^H({\bf{z}})|{\bf{z}}\in\gamma_n\}$. Then we  pick the point (from the candidate set $\Gamma_N$) so that  the  corresponding low-fidelity solution is farthest away from the existing spanned low-fidelity approximation space $U^L(\gamma_k)$, to be the next sampling point: 
\begin{equation}\label{greedy}
{\bf{z}}_{i_{k+1}} = \arg \max_{{\bf{z}}\in \Gamma_N} d^L(u^L({\bf{z}}), U^L(\gamma_k)), \quad \gamma_{k+1} = \gamma_k \cup {\bf{z}}_{i_{k+1}}, 
\end{equation}
where $d^L(v,W)$ is the distance between a function $v\in u^L(\Gamma_N)$ and the space $W\in u^L(\gamma_k)$. We then repeat this step to select all $n$ important points. The whole procedure can be efficiently implemented by apply the pivoted Cheloskey decompostion on $u^L(\Gamma_N)$^^>\cite{NGX14, ZNX14}. 

\bigskip

{\bf Bi-fidelity approximation}. 
For any given $z$, to efficiently compute bi-fidelity approximation in \eqref{UB_1}, it is desirable to find a cheap yet reasonably accurate approximation of $c_k({\bf z})$. The bi-fidelity approach developed in^^>\cite{NGX14, ZNX14} learn these coefficients  from the inexpensive low-fidelity model $u^{L}({\bf z})$. Specifically,
for any given ${\bf z}$, we shall compute the low-fidelity solution $u^L({\bf{z}})$ and then
construct its best approximation in the spanned low-fidelity approximation
space $U^L(\gamma_n)$ by orthogonal projection:
\begin{equation}\label{C-N}
u^L({\bf{z}}) \approx \mathcal{P}_{U^L(\gamma_n)}u^L({\bf{z}}) = \sum_{k=1}^n c_k^L({\bf{z}})u^L({\bf{z}}_{k}), \quad {\bf{z}}_{k}\in \gamma_n, 
\end{equation}
where the low-fidelity projection coefficients can be computed as follows:
\begin{equation}\label{Gc}
{\bf G}^L {\bf c}^L = {\bf f}, \qquad {\bf f}^L = (f_k^L)_{1\leq k\leq n}, \qquad f_k^L = 
\langle u^L({\bf{z}}), u^L({\bf{z}}_{k})\rangle,
 \end{equation}
and
 ${\bf G}^L$ is the Gramian matrix of $u^L(\gamma_n)$, 
\begin{equation}\label{GM}
 ({\bf G}^L)_{ij} =  \left\langle u^L({\bf{z}}_i), u^L({\bf{z}}_j) \right\rangle^L, \qquad 1 \leq i,\, j \leq n, 
\end{equation}
where  $\langle\cdot,\cdot\rangle^L$ is the standard inner product associated with $U^L(\gamma_n)$. 

Under certain condition^^>\cite{NGX14, ZNX14}, $c^L_k({\bf{z}})$ can be reasonably good approximation for the high-fidelity coefficients in \eqref{UB_1}. Consequently,  we can  construct the bi-fidelity approximation of the high-fidelity approximation
solution $u^H({\bf{z}})$ as follows: 
\begin{equation}\label{UB}
 u^B({\bf{z}}) = \sum_{k=1}^nc_k^L({\bf{z}}) u^H({\bf{z}}_{k}),
\end{equation}

To put things together, we outline the major steps for the bi-fidelity approximation of the high-fidelity sample for a given ${\bf{z}}$  in {Algorithm~\ref{BiFi-pod}}. 

\begin{remark}
As we mentioned above, the bi-fidelity method relies on the assumption that the low-fidelity coefficients are similar to the high-fidelity coefficients in the parameter space under certain conditions, as stated in^^>\cite{NGX14, ZNX14}. This may not hold for some problems. In this case, a correction mapping between low-fidelity and high-fidelity coefficient can be constructed by leveraging approximation power of neural network to further improve the accuracy, if the additional high-fidelity data are available.  We refer readers to^^>\cite{LZ21} for details on this approach.
\end{remark}

To construct the bi-fidelity approximation of the high-fidelity mean, with the bi-fidelity surrogate $u^B({\bf z})$, we can employ the Monte Carlo or other quadrature-based methods to compute the statistical moments quickly by evaluating bi-fidelity surrogates. 
To further reduce the number of bi-fidelity surrogate evaluations, 
a more efficient extension to this approach is developed in^^>\cite{zhu2017multi}. The general procedure is as follows:
\begin{itemize}
    \item Compute the low-fidelity sample mean (via the Monte Carlo or quadrature rules):
\begin{equation}
\mu^L = \sum_{i=1}^{M} w_i u^L({\bf z}_i),
\end{equation}
where the high-fidelity mean $\mu^H$ is defined similarly.
\item Compute its best approximation  on the low-fidelity approximation space $U^L(\gamma_n)$ by orthogonal projection, 
\begin{equation}
\mu^L \approx \mathcal{P}_{U^L(\gamma_n)} \mu^L =\sum_{i=1}^{n} c_i^L u^L({\bf z}_i),
\end{equation}
where the expansion coefficients $c^L$ are computed by solving the following linear system:
\begin{equation}
 {\bf G}^L c^L = g^L, \quad g^L = \left\langle \mu^L, u^L({\bf z}_{j})  \right\rangle ^L , \qquad 1 \leq j \leq n, \, {\bf z}_j \in \gamma_n.
\end{equation}
\item Using this coefficient $c^L_{k}$ as the surrogate of the high-fidelity coefficient of the high-fidelity mean $\mu^H$, the bi-fidelity approximation of the high-fidelity mean can be constructed as follows:
\begin{equation}
\mu^B = \sum_{k=1}^{n} c_{k}^L u^H({\bf z}_{k}),\quad z_k\in \gamma_n. 
\end{equation}
\end{itemize}
Note that in this way, only one bi-fidelity surrogate evaluation is required. The standard deviation can be computed similarly. We refer readers to^^>\cite{zhu2017multi} for additional details.

\begin{algorithm*}[ht!]
\caption{A bi-fidelity approximation for a high-fidelity solution at given $z$}
\label{BiFi-pod}

Given a candidate sample set $\Gamma_N = \{{\bf{z}}_1, {\bf{z}}_2, \hdots, {\bf{z}}_N\}\subset I_z $, run the low-fidelity model $u^L({\bf{z}}_j)$ for each ${\bf{z}}_j \in \Gamma_N$.

Select $n$ ``important'' points  $\gamma_N$ from $\Gamma_N$, where $\gamma_N=\{{\bf{z}}_{i_1}, \cdots {\bf{z}}_{i_n} \} \subset\Gamma_N$ and the low-fidelity approximation space by $U^L(\gamma_n)$.

Run high-fidelity simulation only at the point in the selected sample set $\gamma_n$. 

For any given ${\bf{z}}$, run the low-fidelity solver to get the low-fidelity solution $u^L({\bf{z}})$ and compute its low-fidelity coefficients 
in \eqref{C-N}:
\begin{equation*}
u^L({\bf{z}}) \approx \mathcal{P}_{U^L(\gamma_n)}u^L({\bf{z}}) = \sum_{k=1}^n c_k^L({\bf{z}})u^L({\bf{z}}_k), \quad {\bf{z}}_k \in \gamma_n. \end{equation*}

Construct the bi-fidelity approximation by applying the same approximation rule as in low-fidelity model with \eqref{UB}:
\begin{equation*}
u^B({\bf{z}})  = \sum_{k=1}^n c_k^L({\bf{z}})u^H({\bf{z}}_k). \end{equation*}
\end{algorithm*}

\section{Numerical examples}
To examine the performance of the proposed methodology, two benchmark tests are considered: 
the first concerning the kinetic transport model with the SIR compartmentalization discussed in Section \ref{section_SIR}, and the second one regarding the extension to the SEIAR modeling presented in Section \ref{section_SEIAR}.
It is worth to highlight that, in these tests, we collect all the quantities of interests together in a single vector to choose the ${\bf{z}}_{k}$ points of the bi-fidelity algorithm, for either the 3 compartments $S,I,R$ in Test 1 or the five compartments $S$, $E$, $I$, $A$, $R$ in Test 2. From the experiments, we also find that the results are similar if we choose the ${\bf{z}}_{k}$ points separately for the different epidemic compartments.

\subsection{Test 1: SIR model with heterogeneous environment}
We first consider the initial distributions of the high-fidelity kinetic SIR model \eqref{eqn:kinetic} as follows:
\begin{equation}
f_i(x,v,0) = c\, i(x,0)\, e^{-\frac{v^2}{2}}, \qquad i \in \{S,I,R\}
\label{eq.ICf}
\end{equation}
where $c=\frac{1}{2} \sum_{i=1}^{N_G} w_i \, e^{-\frac{\zeta_i^2}{2}}$ is a re-normalization constant, with $N_G$ number of Gauss-Legendre quadrature points as defined in Section \ref{sect_AP_IMEX}, and
$$ S(x,0)= 1 - I(x,0), \qquad I(x,0) = 0.01 e^{-(x-10)^2}, \qquad R(x,0)=0, $$
in the physical domain $L=[0,20]$. The initial fluxes $J_S(x,0), J_I(x,0)$ and $J_R(x,0)$ are null and we consider periodic boundary conditions. The same initial conditions for $S, I, R$ and $J_S, J_I, J_R$ are imposed in the low-fidelity SIR model \eqref{eq:density}-\eqref{eq.SIR_kinetic-fluxes}. In the following tests, we consider a 2-dimensional random vector ${\bf{z}}=(z_1, z_2)^T$, with independent random parameters $z_1$ and $z_2$ that follow a uniform distribution, $z_j \sim \mathcal{U}(-1,1)$, $j=1,2$.

We analyze the behavior of the proposed methodology concerning spatially heterogeneous environments, considering a spatially variable contact rate^^>\cite{Bert,Wang2020} 
\begin{equation*}
\beta(x,{\bf{z}}) = \beta^0({\bf{z}}) \left(1+0.05\sin\left(\frac{13\pi x}{20}\right)\right),
\end{equation*}
where 
$$\beta^0({\bf{z}}) = 11(1 + 0.6 z_1),$$ and a recovery rate 
$$\gamma({\bf{z}}) = 10(1 + 0.4 z_2).$$
In the incidence function, we set $\kappa=0$ and $p=1$. Therefore, we simulate an infectious disease characterized by an $R_0$ perturbed around the value 1 depending on the random fluctuations. We discretize the spatial domain, both in the high-fidelity and low-fidelity models, with $N_x=150$ cells and consider $N_v=8$ velocities for the high-fidelity model. To compute the reference solutions for the mean and standard deviation, we use a $3$-rd level sparse grid quadrature based on Clenshaw-Curtis rules for the choice of the stochastic collocation nodes, with a total of $29$ points, for both the high-fidelity and low-fidelity models. 

\paragraph{Test 1 (a):}
In this case, a parabolic configuration of speeds and relaxation parameters is considered, setting $\lambda_i^2=10^5$, $i \in \{S,I,R\}$, and
$\tau_i=10^{-5}$ in the low-fidelity model and $\tau_i=3\times 10^{-5}$ in the high-fidelity model, to maintain consistency of the two simulations (see Remark~\ref{remark_consistency}). The time step size results $\Delta t = 0.89 \times 10^{-2}$ in both models, nevertheless the low-fidelity code is almost 5 times faster than the high-fidelity one ($t_{CPU}^{HF} \approx 70.0 \,s$, $t_{CPU}^{LF} \approx 14.3 \,s$ to simulate until $t=5$).

\begin{figure}[!t]
\centering
\includegraphics[scale=0.38]{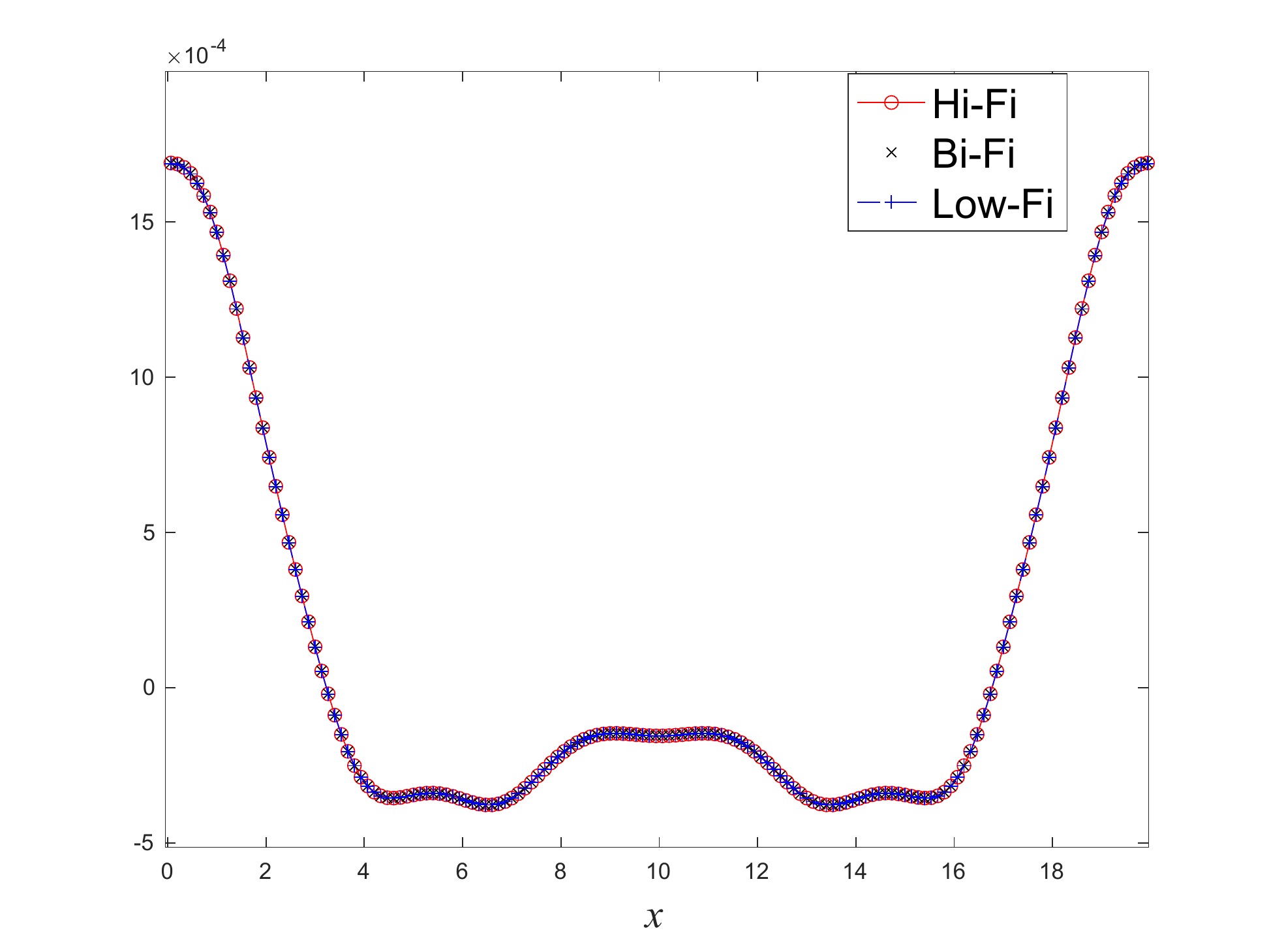}
\includegraphics[scale=0.38]{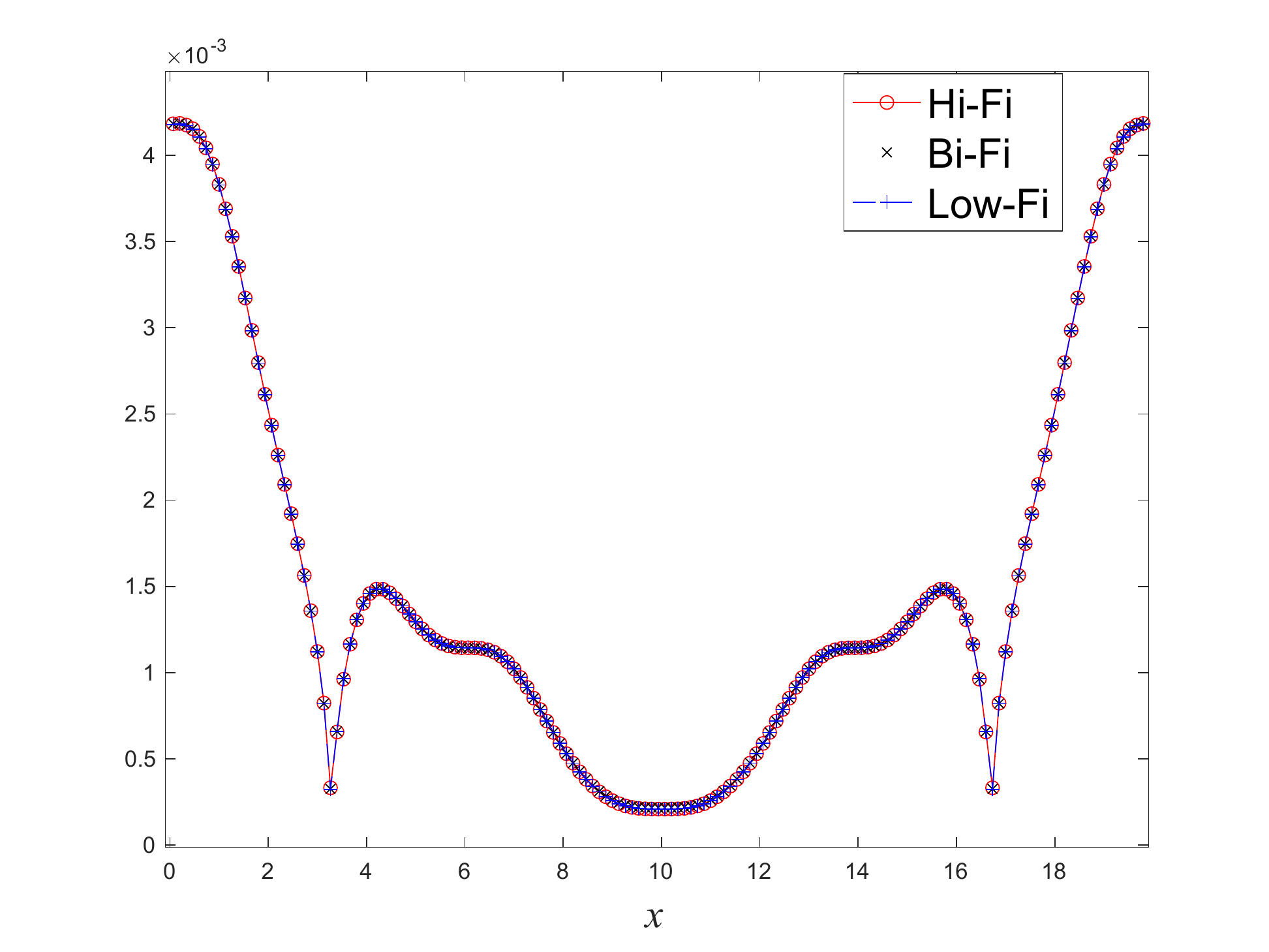}
\includegraphics[scale=0.38]{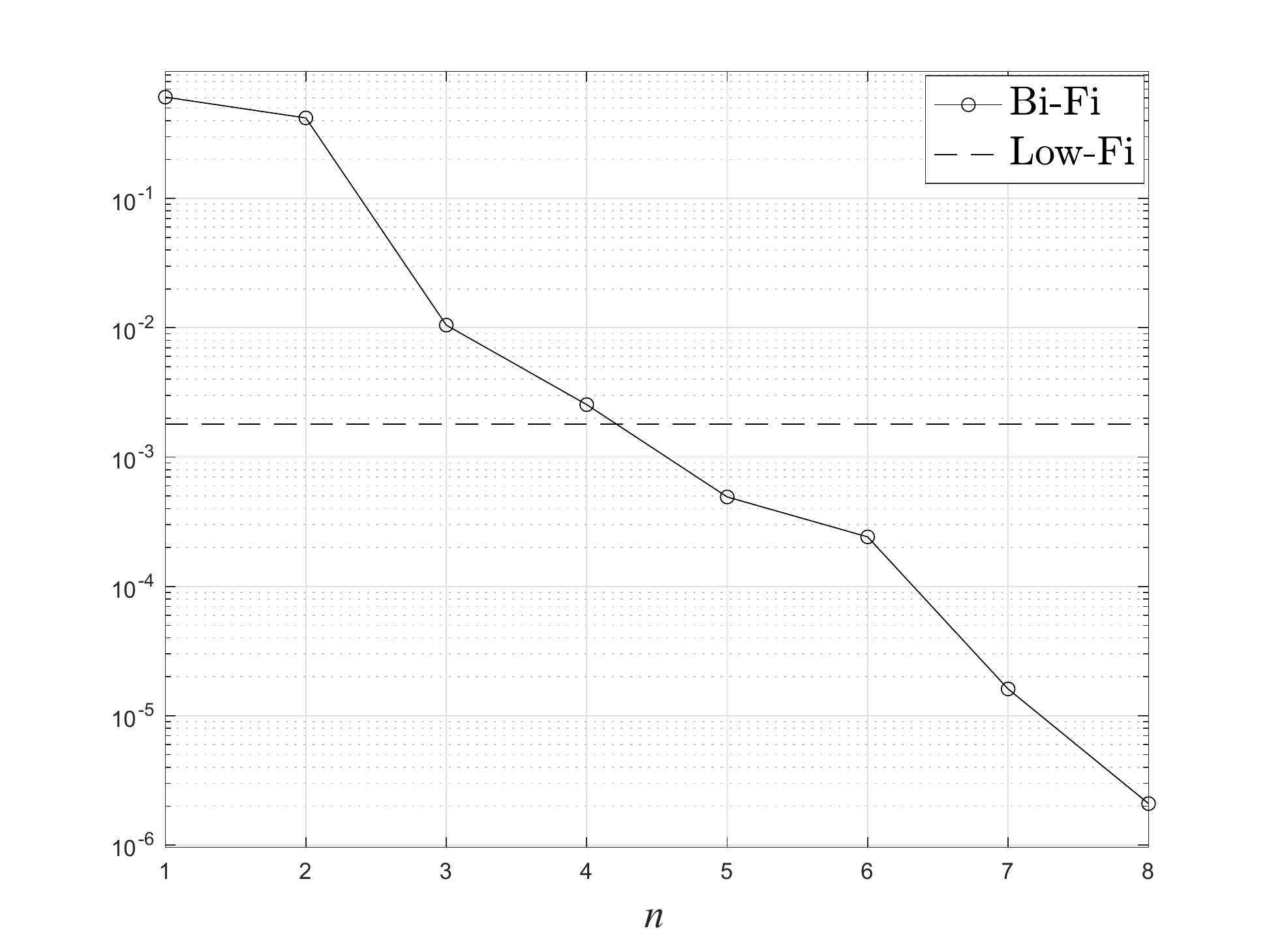}
\includegraphics[scale=0.38]{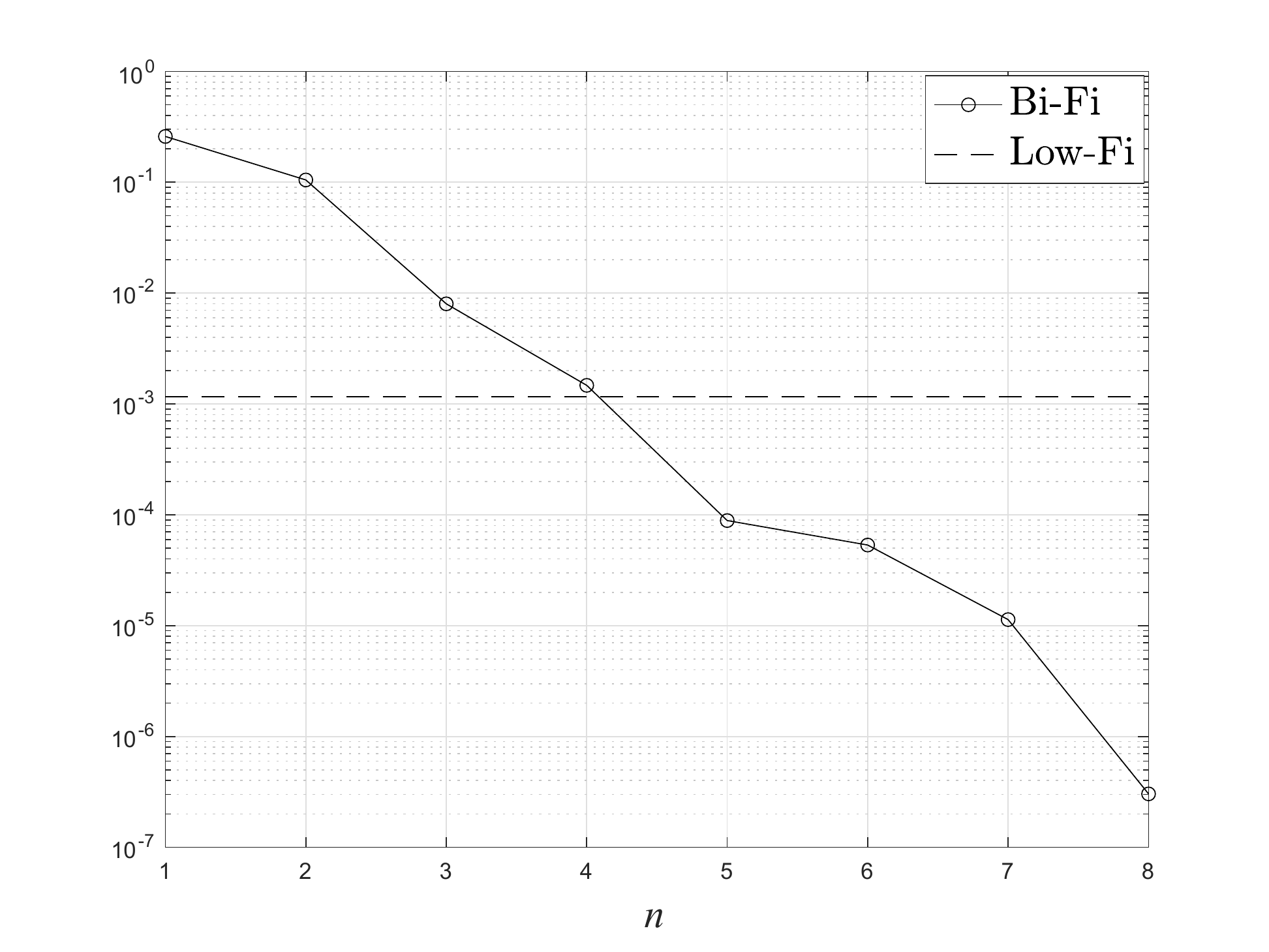}
\caption{Test 1 (a): SIR model in diffusive regime. First row: expectation (left) and standard deviation (right) obtained at $t=5$ for the variable $I$ with the three methodologies, by using $n=8$ points for the bi-fidelity approximation. Second row: relative $L^2$ errors of the bi-fidelity approximation for the mean (left) and standard deviation (right) of density $I$ with respect to the number of "important" points $n$ used in the bi-fidelity algorithm, compared with low-fidelity errors.}
\label{fig.test1a}
\end{figure}

In the first row of Figure \ref{fig.test1a}, the expectation and standard deviation of the solution of compartment $I$ for the high-fidelity model, the low-fidelity model and  the bi-fidelity approximation at time $t=5$ are shown. 
As expected, since high-fidelity and low-fidelity models share the same diffusive limit, a perfect agreement of the solutions can be observed. To confirm this, in the second row of the same figure, we also plot the $L^2$ errors of the mean and standard deviation between the bi-fidelity and high-fidelity solutions at $T=5$ with respect to the number of selected "important" points $n$ of the bi-fidelity algorithm. A fast error decay is clearly observed.
With only $n=8$ hi-fidelity sample points, bi-fidelity approximation can achieved a relative error of $\mathcal{O}(10^{-6})$ for both the mean and standard deviation. 

\begin{figure}[!t]
\centering
\includegraphics[scale=0.38]{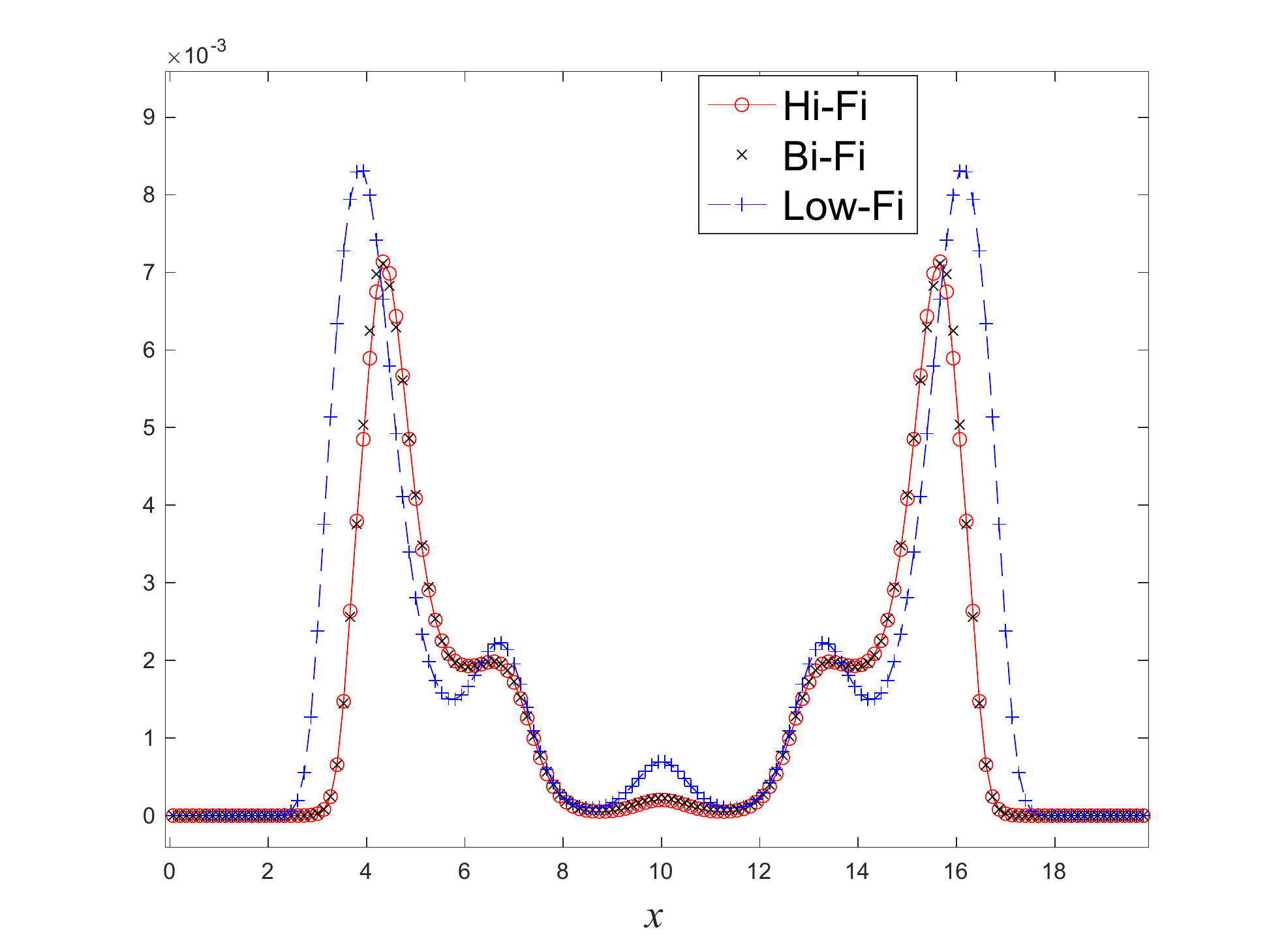}
\includegraphics[scale=0.38]{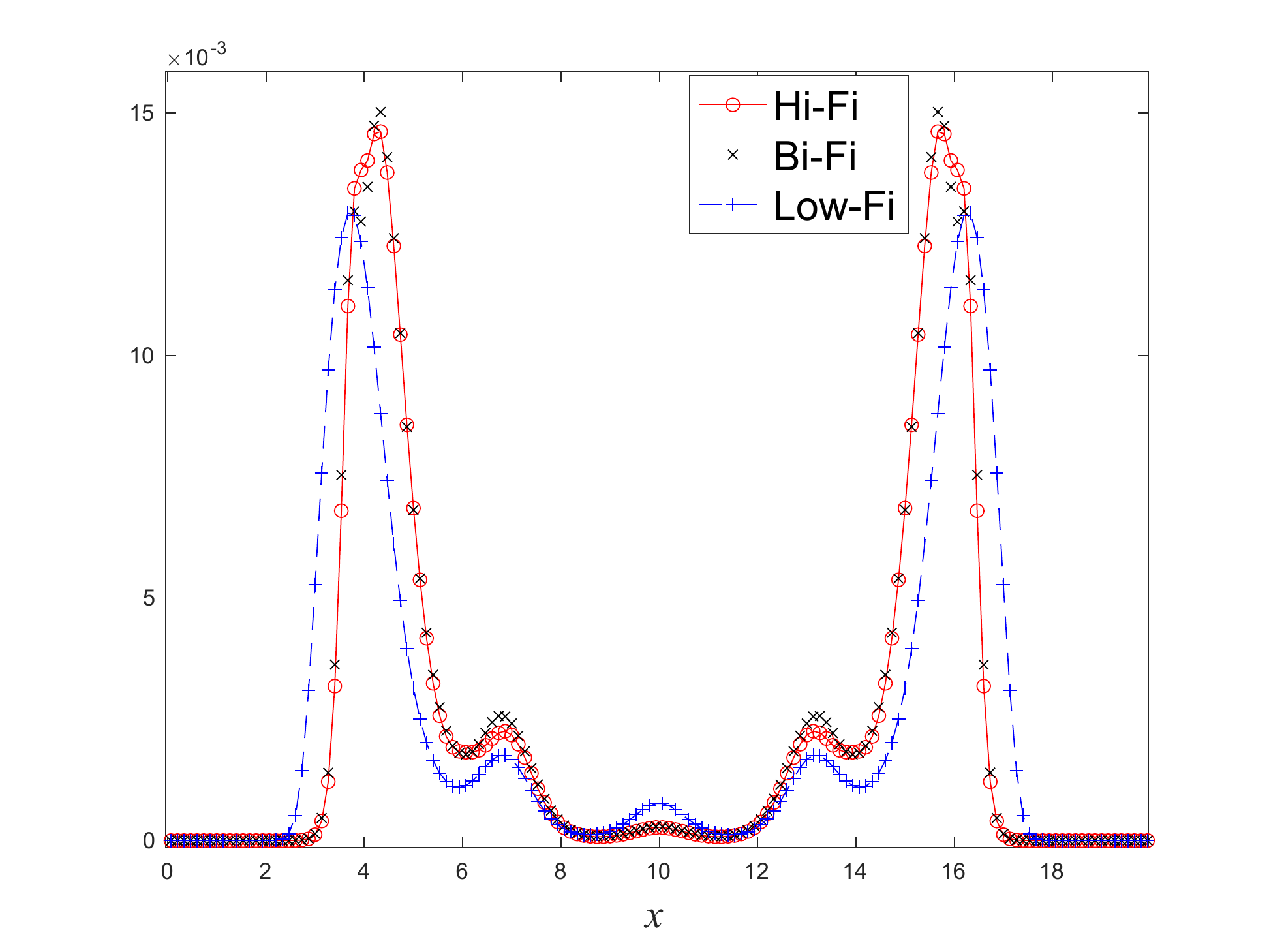}
\includegraphics[scale=0.38]{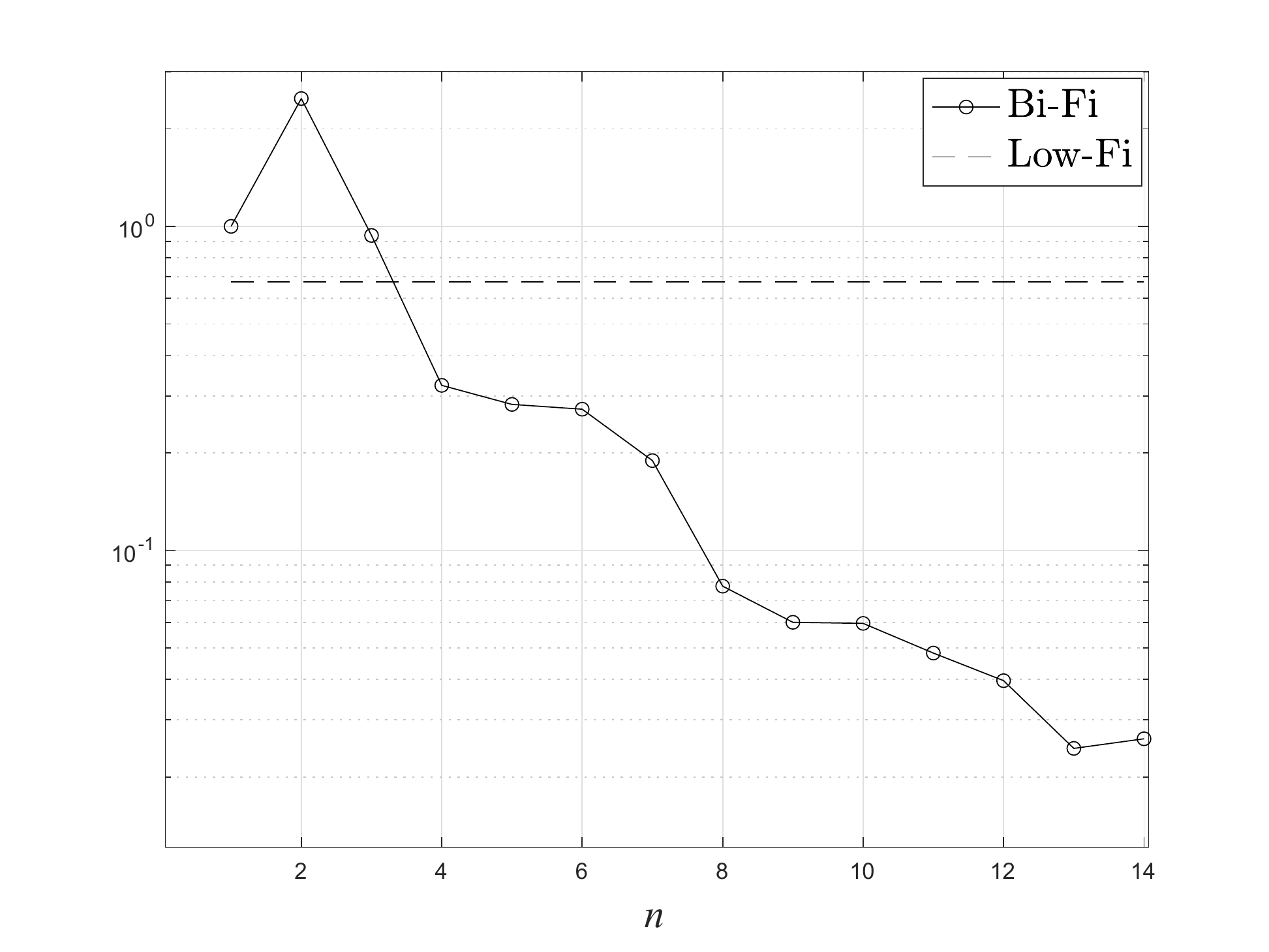}
\includegraphics[scale=0.38]{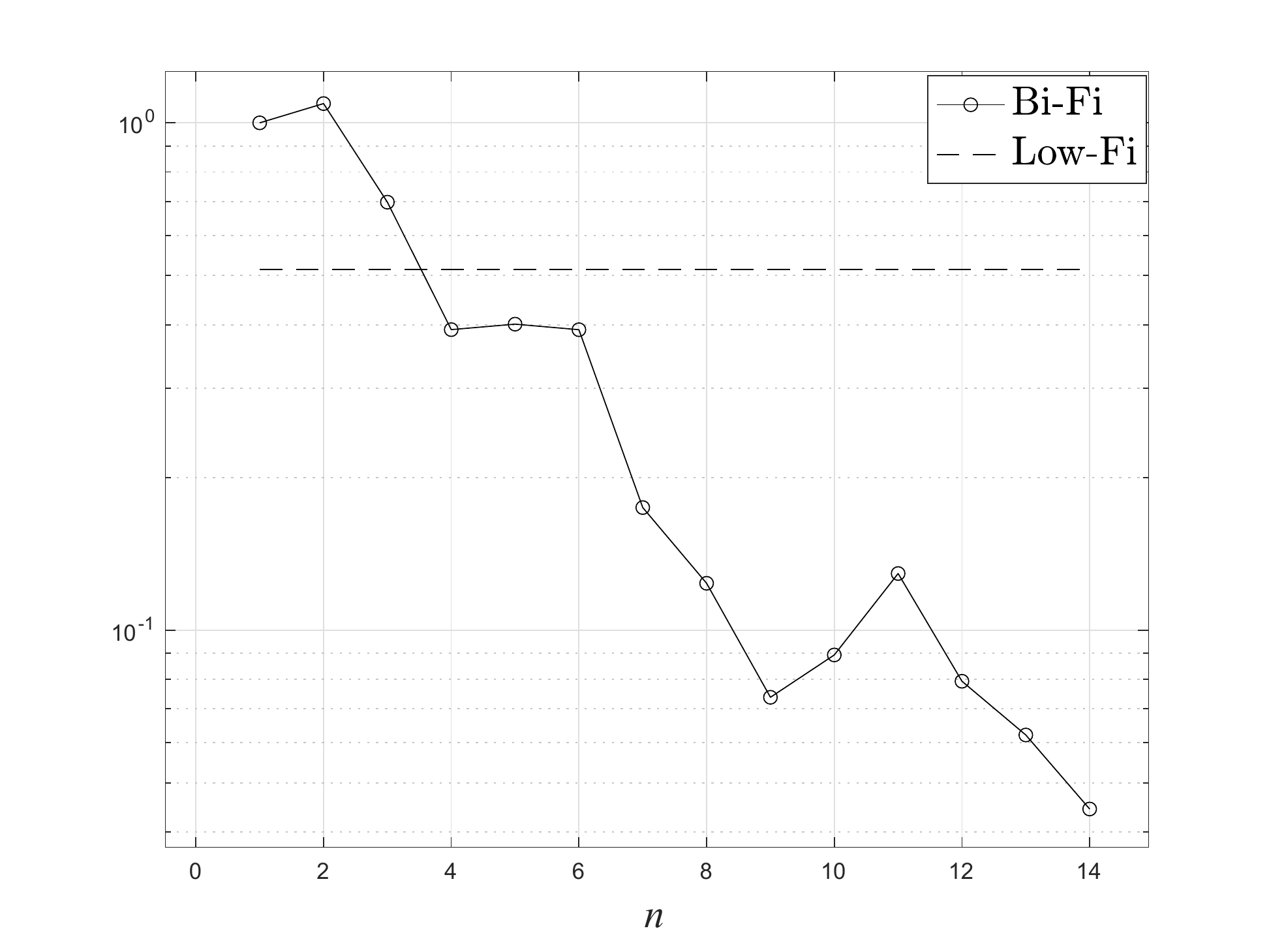}
\caption{Test 1 (b): SIR model in hyperbolic regime. First row: expectation (left) and standard deviation (right) obtained at $t=5$ for the variable $I$ with the three methodologies, by using $n=14$ points for the bi-fidelity approximation. Second row: relative $L^2$ errors of the bi-fidelity approximation for the mean (left) and standard deviation (right) of density $I$ with respect to the number of "important" points $n$ used in the bi-fidelity algorithm, compared with low-fidelity errors.}
\label{fig.test1b}
\end{figure}
\paragraph{Test 1 (b):}
In this second case, we consider the hyperbolic regime by letting $\lambda_i=1$, $i \in \{S,I,R\}$, while $\tau_i=1$ in the low-fidelity model and $\tau_i=3$ in the high-fidelity model, for consistency. The time step size results $\Delta t = 0.12$ in both models, but the low-fidelity simulation is again almost 5 times faster than the high-fidelity one ($t_{CPU}^{HF} \approx 5.6\, s$, $t_{CPU}^{LF} \approx 1.3\, s$ to simulate until $t=5$).

In Figure \ref{fig.test1b}, the results of Test 1(b) are reported for the infectious compartment $I$ at time $t=5$. The first row shows expectation and standard deviation of the high-fidelity, the low-fidelity and  the bi-fidelity solutions. 
While low-fidelity model fails to capture details of the hi-fidelity solutions around peaks due to the different velocity setting considered in the two models,
bi-fidelity approximation are in almost perfect agreement with the high-fidelity ones, confirming the validity of the proposed methodology even with hyperbolic configurations of the scaling parameters. In addition,  we plot the relative $L^2$ error convergence of the mean and standard deviation of the bi-fidelity and high-fidelity solutions at $t=5$ with the number of high-fidelity samples, in the second row of Figure \ref{fig.test1b}. Again, we observe a fast error decay for mean and standard deviation with $\mathcal{O}(10)$ high-fidelity samples.


\begin{figure}[ht!]
\centering
\includegraphics[width=0.48\textwidth]{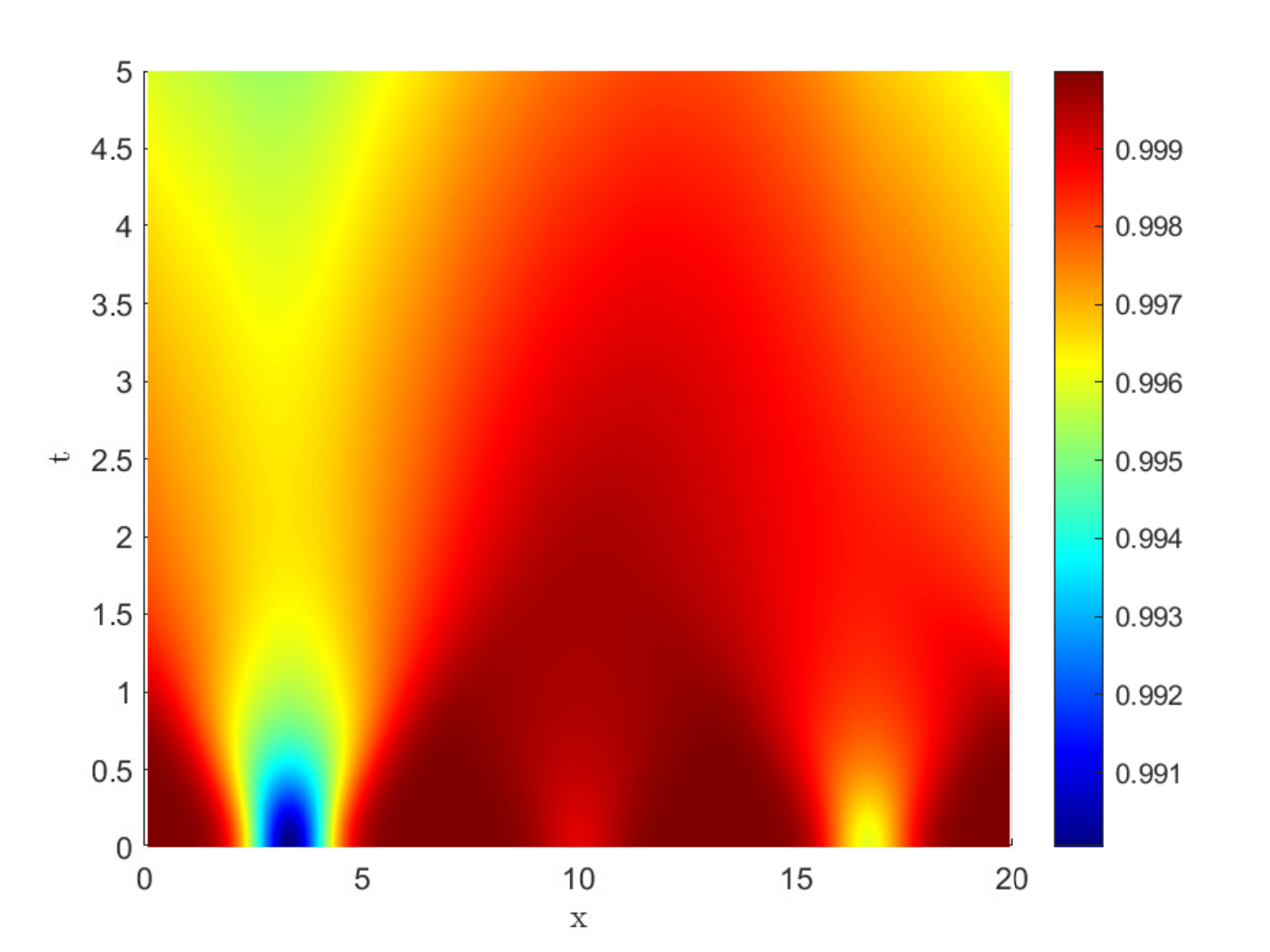}
\includegraphics[width=0.48\textwidth]{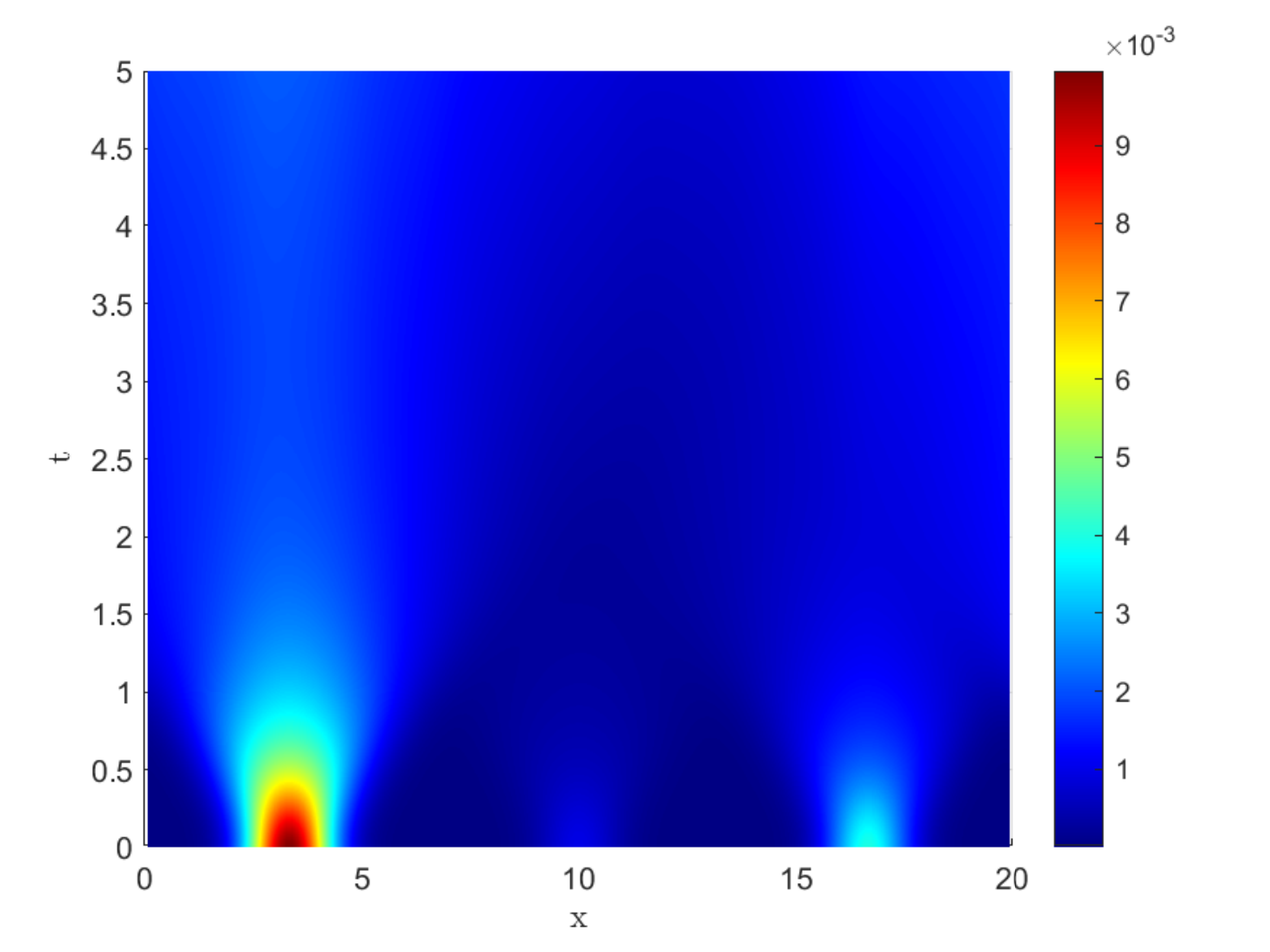}
\includegraphics[width=0.48\textwidth]{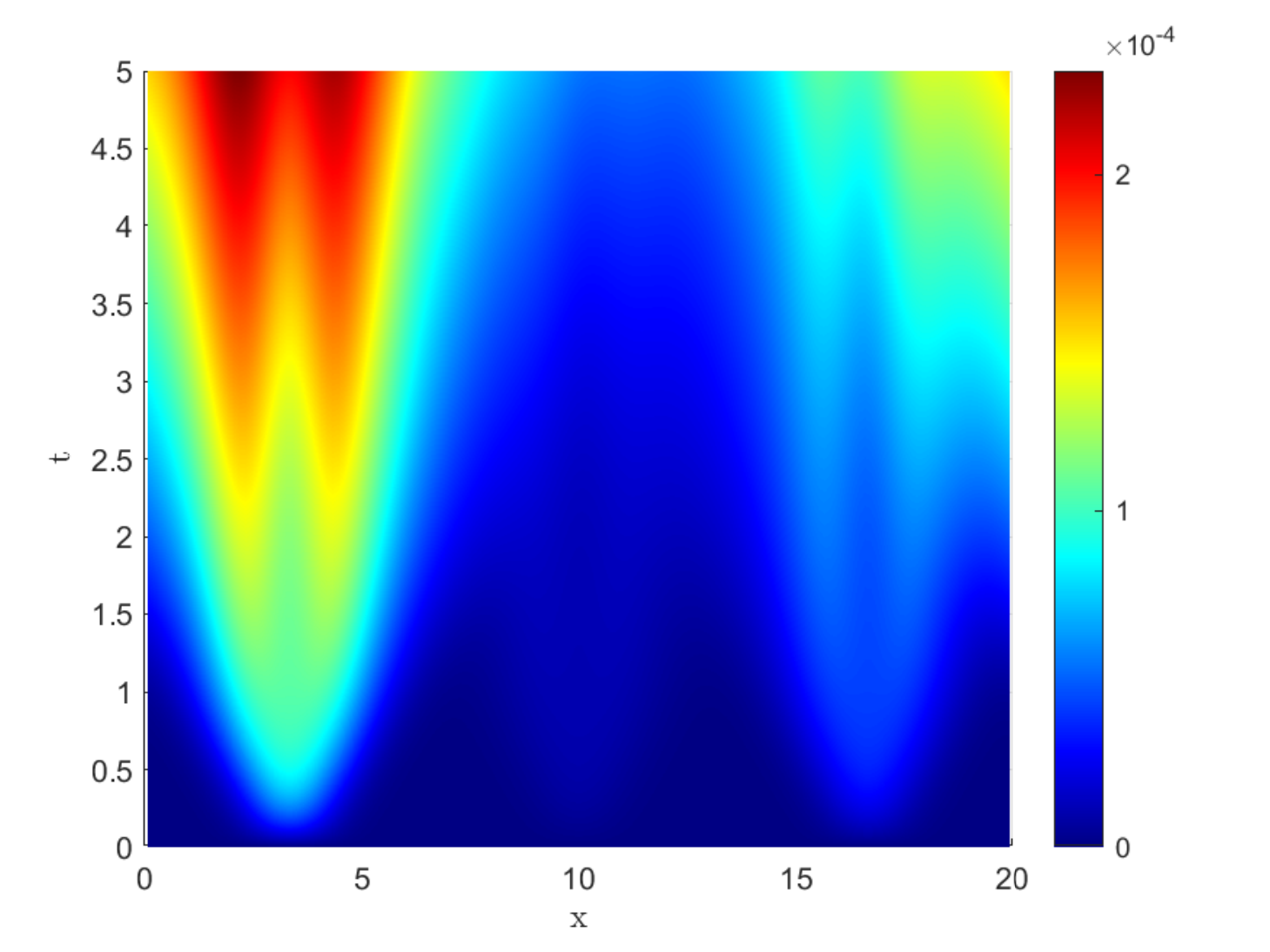}
\includegraphics[width=0.48\textwidth]{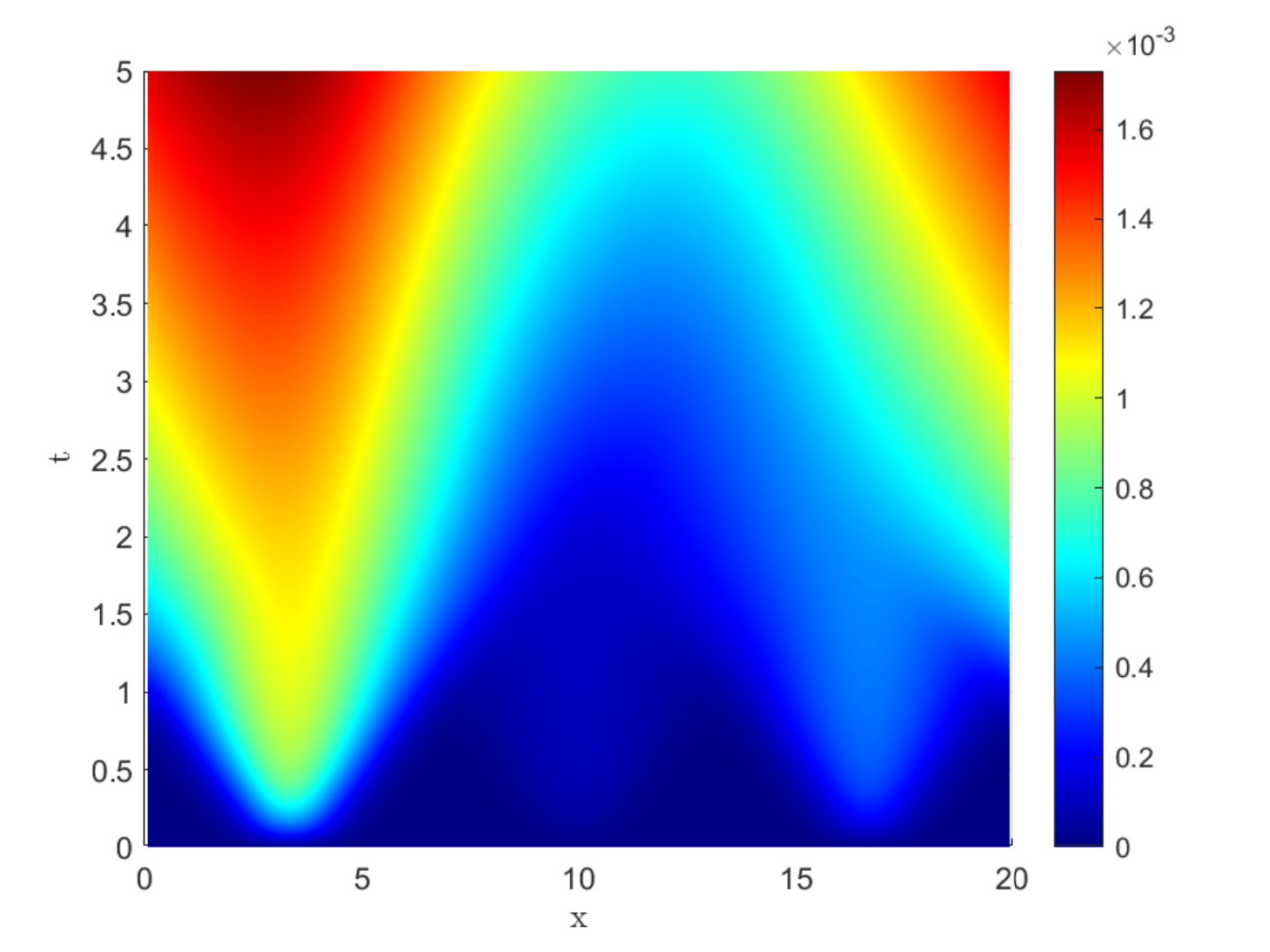}
\caption{Test 2 (a): SEIAR model in intermediate regime. The baseline temporal and spatial evolution of compartments $S$ (first row, left), $E$ (first row, right), $I$ (second row, left) and $A$ (second row, right) in the high-fidelity model.}
\label{fig.test2a_xt}
\end{figure}
\begin{figure}[p]
\centering
\includegraphics[scale=0.38]{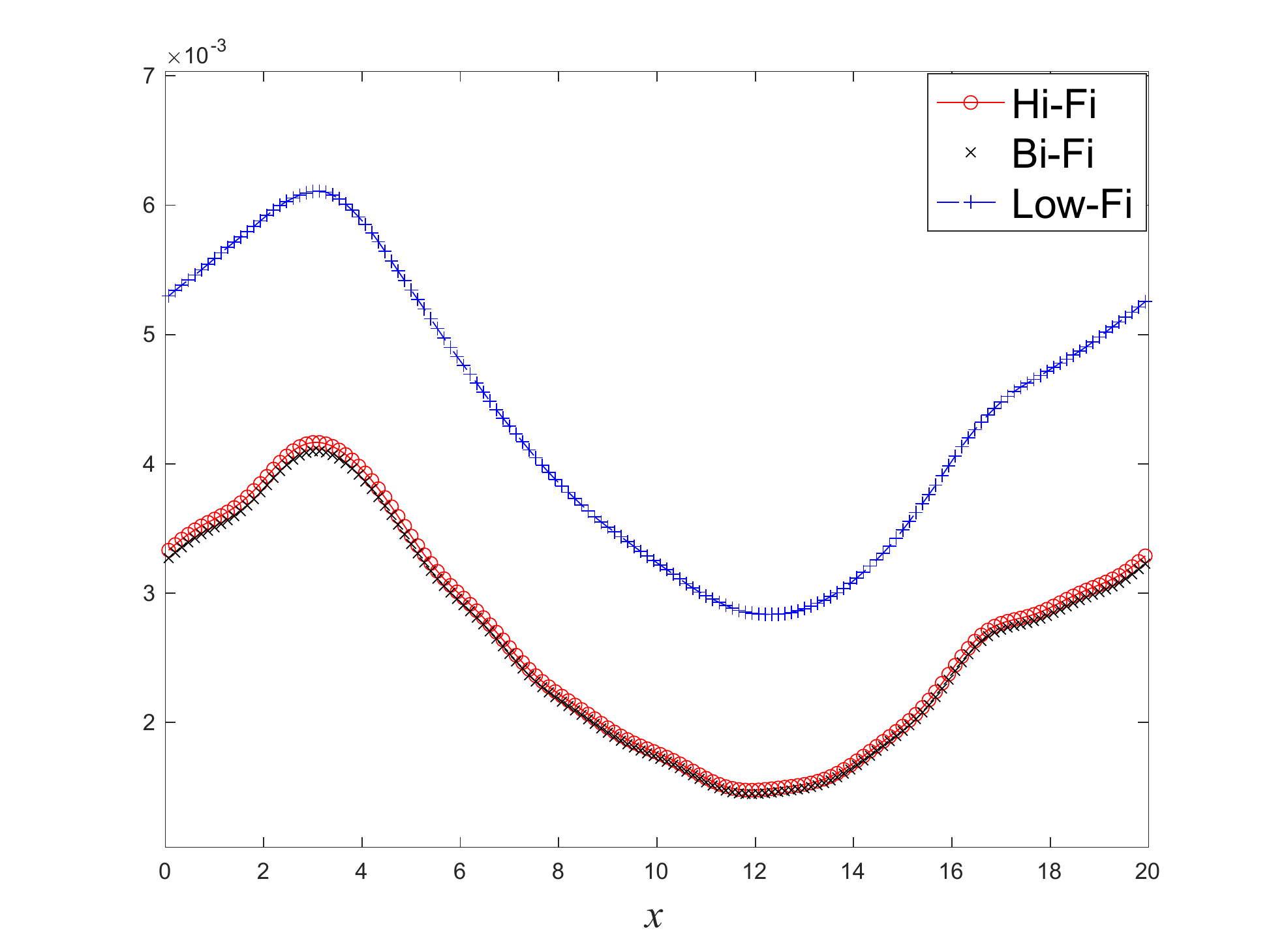}
\includegraphics[scale=0.38]{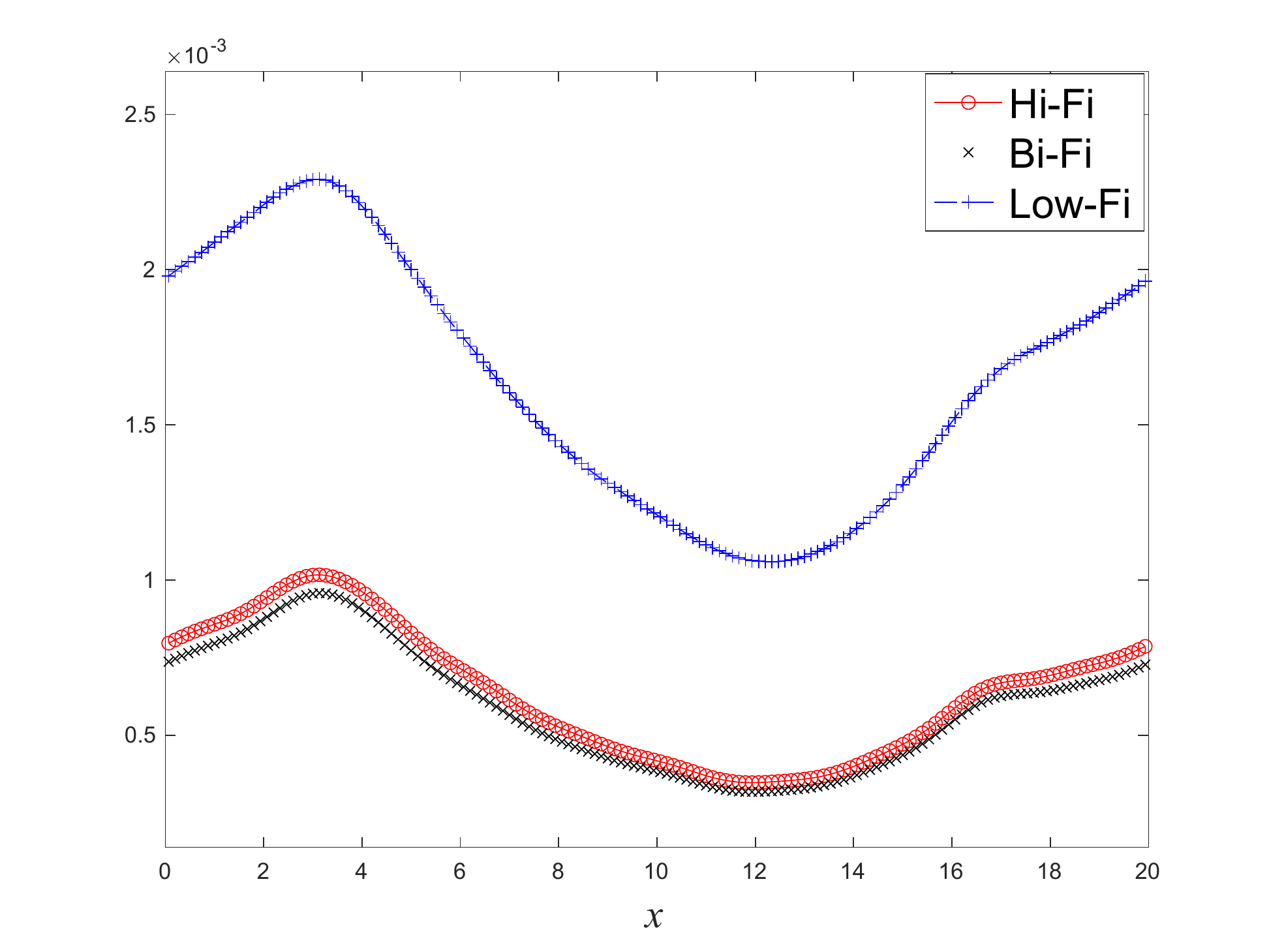}
\includegraphics[scale=0.38]{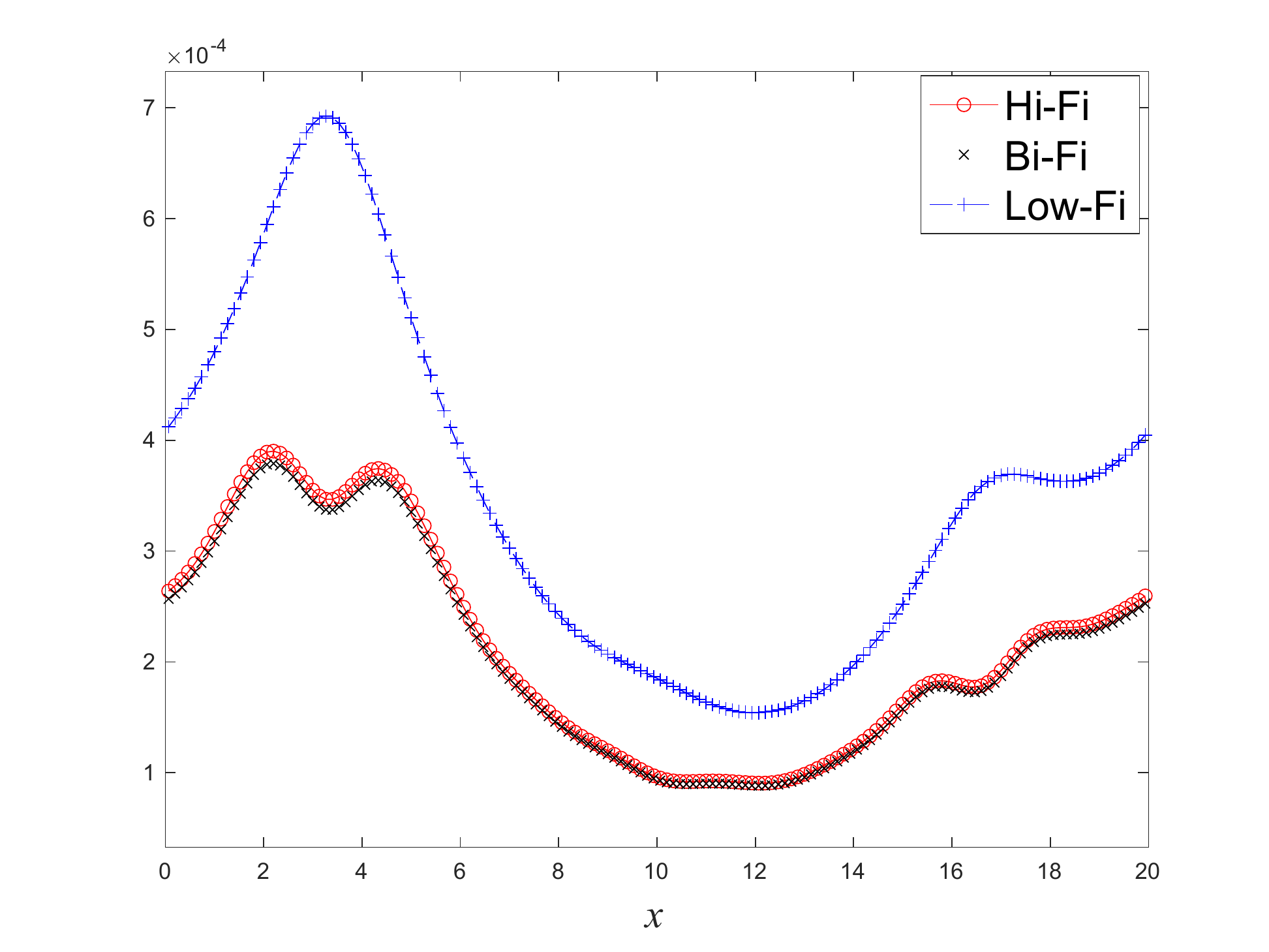}
\includegraphics[scale=0.38]{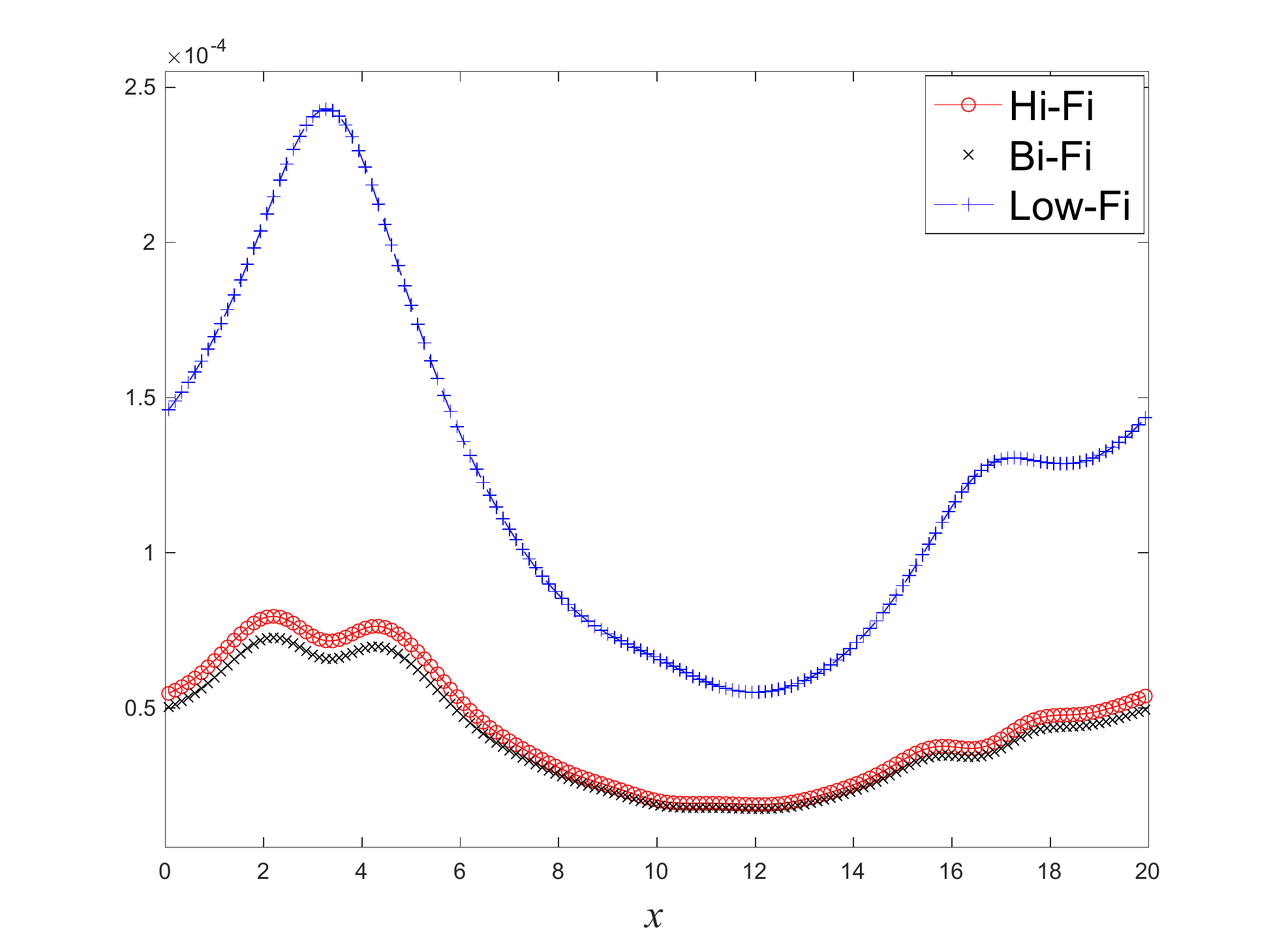}
\includegraphics[scale=0.38]{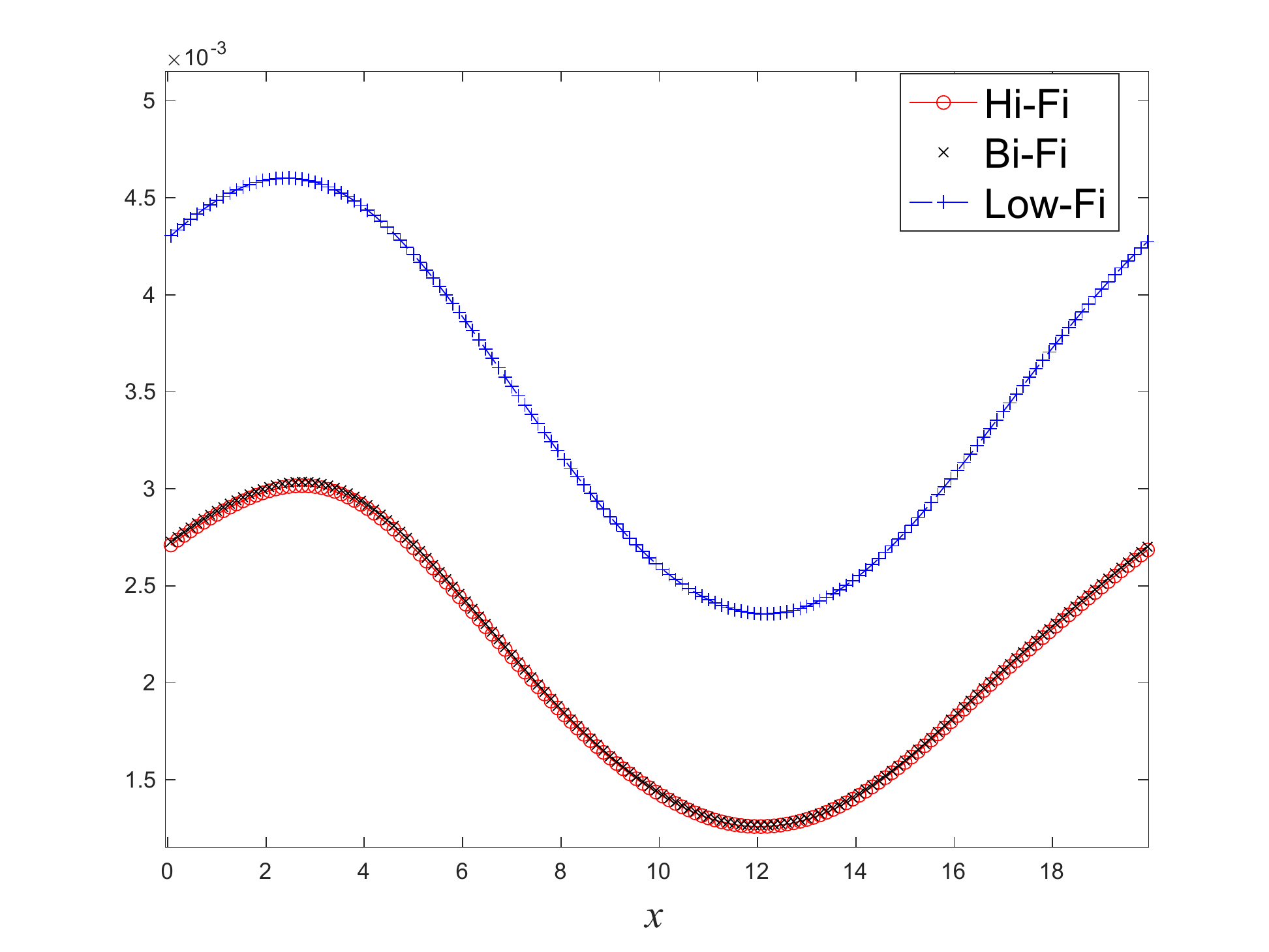}
\includegraphics[scale=0.38]{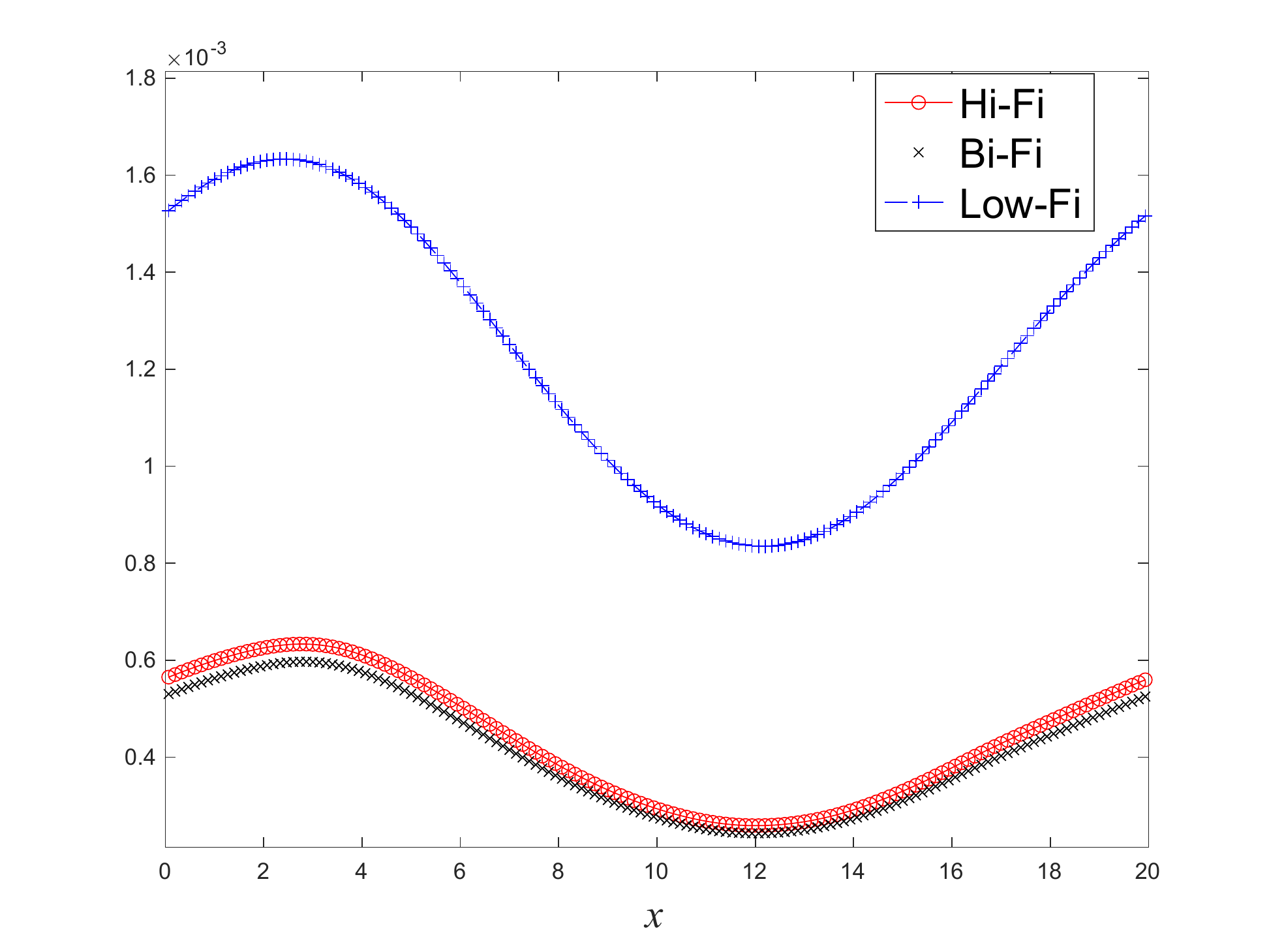}
\caption{Test 2 (a): SEIAR model in intermediate regime. Expectation (left) and standard deviation (right) of densities $E$ (first row), $I$ (second row) and $A$ (third row) at time $t=5$, obtained with the three methodologies, using $n=6$ for the bi-fidelity solution. }
\label{fig.test2a}
\end{figure}
\begin{figure}[ht!]
\centering
\includegraphics[scale=0.38]{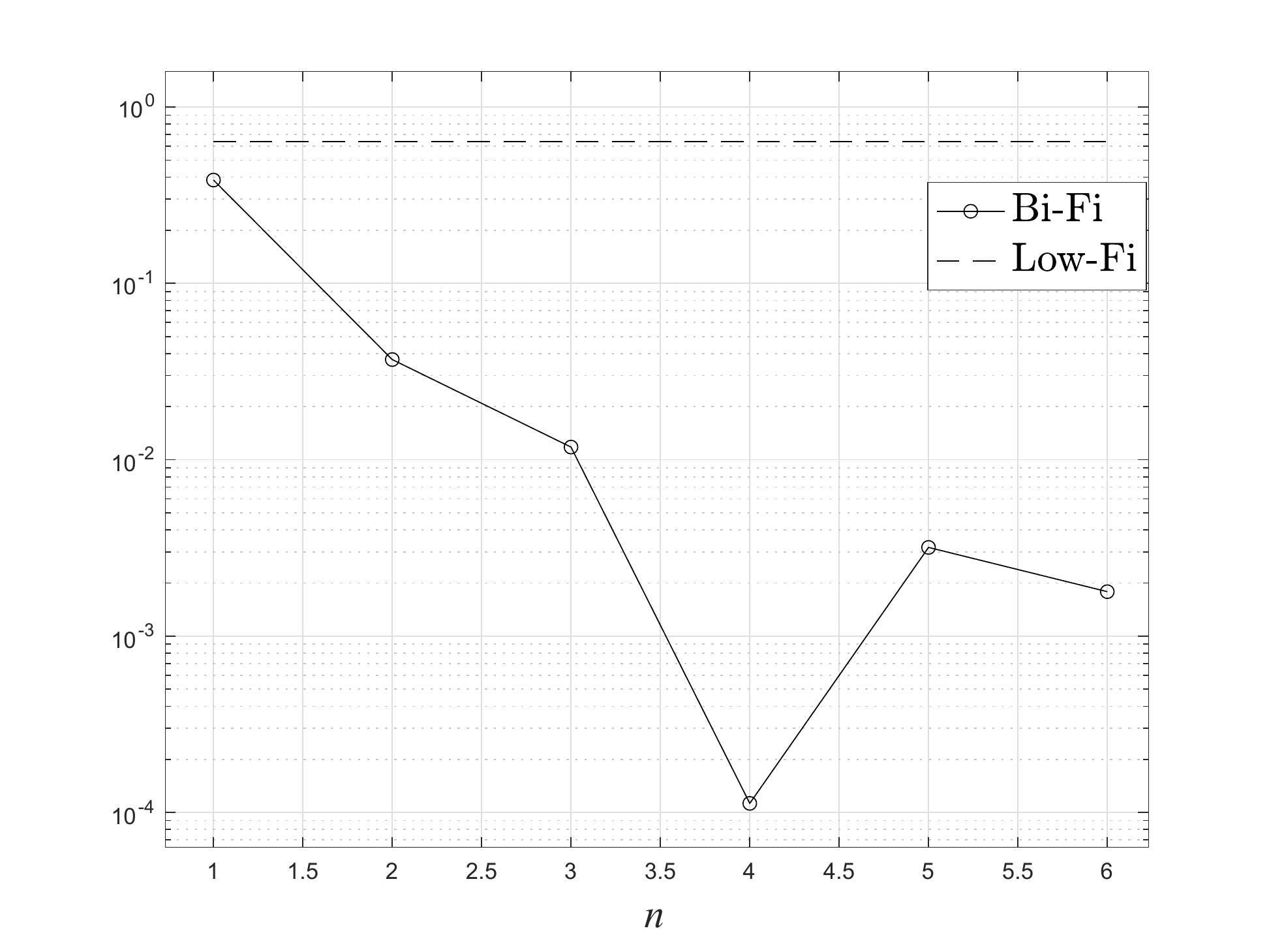}
\includegraphics[scale=0.38]{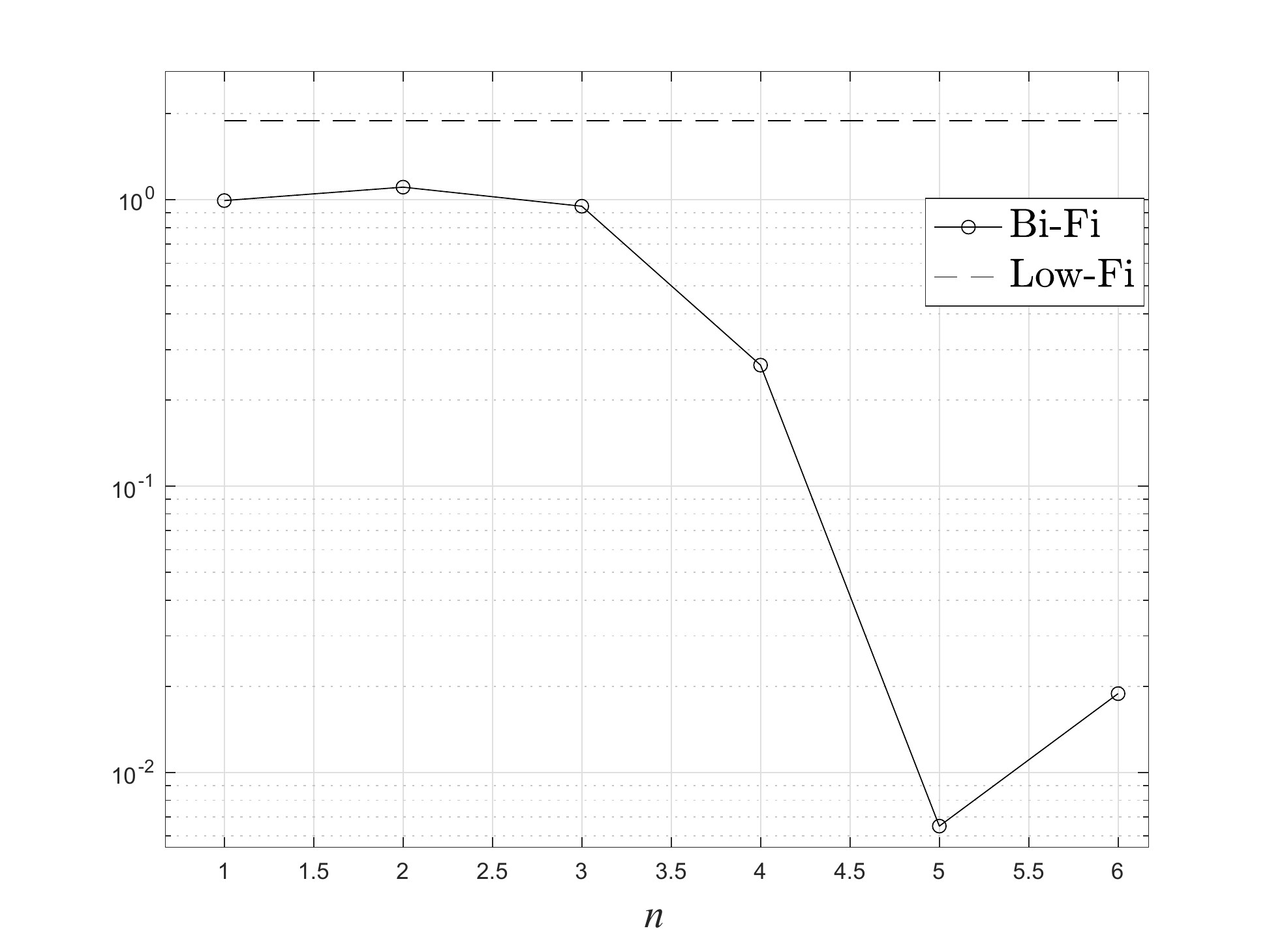}
\caption{Test 2 (a): SEIAR model in intermediate regime. Relative $L^2$ error decay of the bi-fidelity approximation of expectation (left) and standard deviation (right) for the density $A$ with respect to the number of selected "important" points $n$, compared with low-fidelity errors. }
\label{fig.test2a_err}
\end{figure}
\subsection{Test 2: SEIAR model with distinguished epidemic hotspots}
Next, we analyze the effectiveness of the proposed methodology also with the extended SEIAR compartmentalization examining a more realistic epidemic scenario. Let us now consider a 2-dimensional random vector ${\bf{z}}$, this time with $z_j \sim \mathcal{U}(0,1)$, $j=1,2$. We design an initial condition for the low-fidelity SEIAR model \eqref{eq:density-1Dkinetic-SEIAR}-\eqref{eq.SEIARmacro-fluxes} that simulates the presence of 3 cities aligned in the spatial domain $L=[0,20]$ with a different number of exposed (infected but not yet infectious) individuals, subjected to randomness, distributed following a Gaussian function,
$$ E(x,0,{\bf z})= \alpha_1({\bf z})\,e^{-(x-x_1)^2} + \alpha_2({\bf z})\,e^{-(x-x_2)^2} + \alpha_3({\bf z})\,e^{-(x-x_3)^2},$$
where $x_1=10/3$, $x_2=10$, $x_3=50/3$ are the coordinates of the city centers and 
$\alpha_1=0.01(1+z_1)$, $\alpha_2=0.001(1+z_1)$, $\alpha_3=0.004(1+z_1)$ are the stochastic amplitudes. 
In this case, we set the large stochastic amplitude to model large uncertainty present on the initial amount of exposed people when an epidemic starts spreading. Indeed, it is worth to underline that the initial number of exposed individuals is certainly one of the variables majorly affected by uncertainty during epidemic outbreaks^^>\cite{Gatto}. We consider that there are no infectious people at $t=0$, being
$$ S(x,0,{\bf z})=1 - E(x,0,{\bf z}), \qquad I(x,0)=0, \qquad A(x,0)=0, \qquad R(x,0)=0, $$
and $J_S(x,0)=J_E(x,0)=J_I(x,0)=J_A(x,0)=J_R(x,0)=0$, with periodic boundary conditions, to allow a connection also between cities 1 and 3.
The initial distributions of the high-fidelity kinetic SEIAR model \eqref{eq:kineticc-SEIAR} then read as defined in \eqref{eq.ICf}, for $i \in \{S,E,I,A,R\}$.

Concerning epidemic parameters, to simulate the more challenging scenario in which the incidence function presents sinusoidal oscillations in space as well as being greater in the most populated areas, we consider the following distribution of the contact rate related to asymptomatic individuals:
\begin{equation*}
\beta_A(x,{\bf{z}}) = \beta_A^0({\bf{z}}) \left(1 + \frac{1}{2}\,e^{-(x-x_1)^2} + \frac{1}{4}\,e^{-(x-x_2)^2} + \frac{1}{2}\,e^{-(x-x_3)^2} \right) + 0.05\,\mathrm{sin}(2\pi x),
\end{equation*}
where $\beta_A^0$ is affected by the following random fluctuations:
$$\beta_A^0({\bf{z}}) = 0.5(1 + 0.5 z_2).$$
Assuming that highly infectious subjects are mostly detected in the most optimistic scenario, being subsequently quarantined or hospitalized, we set $\beta_I(x,{\bf{z}}) = 0.03\, \beta_A(x,{\bf{z}})$. Then, we fix $\gamma_I = 1/14$, $\gamma_A=1/7$, $a = 1/3$, $\sigma = 1/12.5$, considering these clinical parameters deterministic according to values adopted in^^>\cite{Bert2,Bert3,Buonomo,Gatto} to simulate the COVID-19 spread. Finally, in the incidence functions, we set $\kappa_I=\kappa_A=0$, hence assuming that initially individuals are not aware of the epidemic outbreak, and $p=1$ to work with the standard bi-linear case. The resulting setting permits to simulate and epidemic characterized by a baseline reproduction number $R_0 = 3$.
We consider $N_x=150$ cells for both high-fidelity and low-fidelity simulations, and use $N_v=8$ velocities for the high-fidelity solution.  

\paragraph{Test 2 (a):}
In the first case, an intermediate regime between parabolic and hyperbolic is considered, setting $\lambda_i^2=10$, $i~\in~\{S,E,A,R\}$, $\lambda_I = 0$ and $\tau_i=0.25$, $i~\in~\{S,E,I,A,R\}$, in the low-fidelity model and $\tau_i=0.75$ in the high-fidelity model. The characteristic speed of compartment $I$ is fixed to zero because we assumed that infectious people with severe symptoms are generally detected and isolated^^>\cite{Bert2,Bert3}. 

The baseline space-time evolution of the compartments is shown in Figure \ref{fig.test2a_xt}, where the fast propagation of the epidemic can be observed, especially starting from the first city on the left, as expected. In Figure \ref{fig.test2a}, the expectation and standard deviation of infected individuals, divided by each compartment, are presented for each methodology adopted, highlighting the validity of the bi-fidelity approach. In Figure \ref{fig.test2a_err}, $L^2$ error decay of the bi-fidelity approximations, with respect to the number of selected "important" points $n$, is compared with the error that would be obtained simply using the low-fidelity model.

\begin{figure}[t!]
\centering
\includegraphics[width=0.48\textwidth]{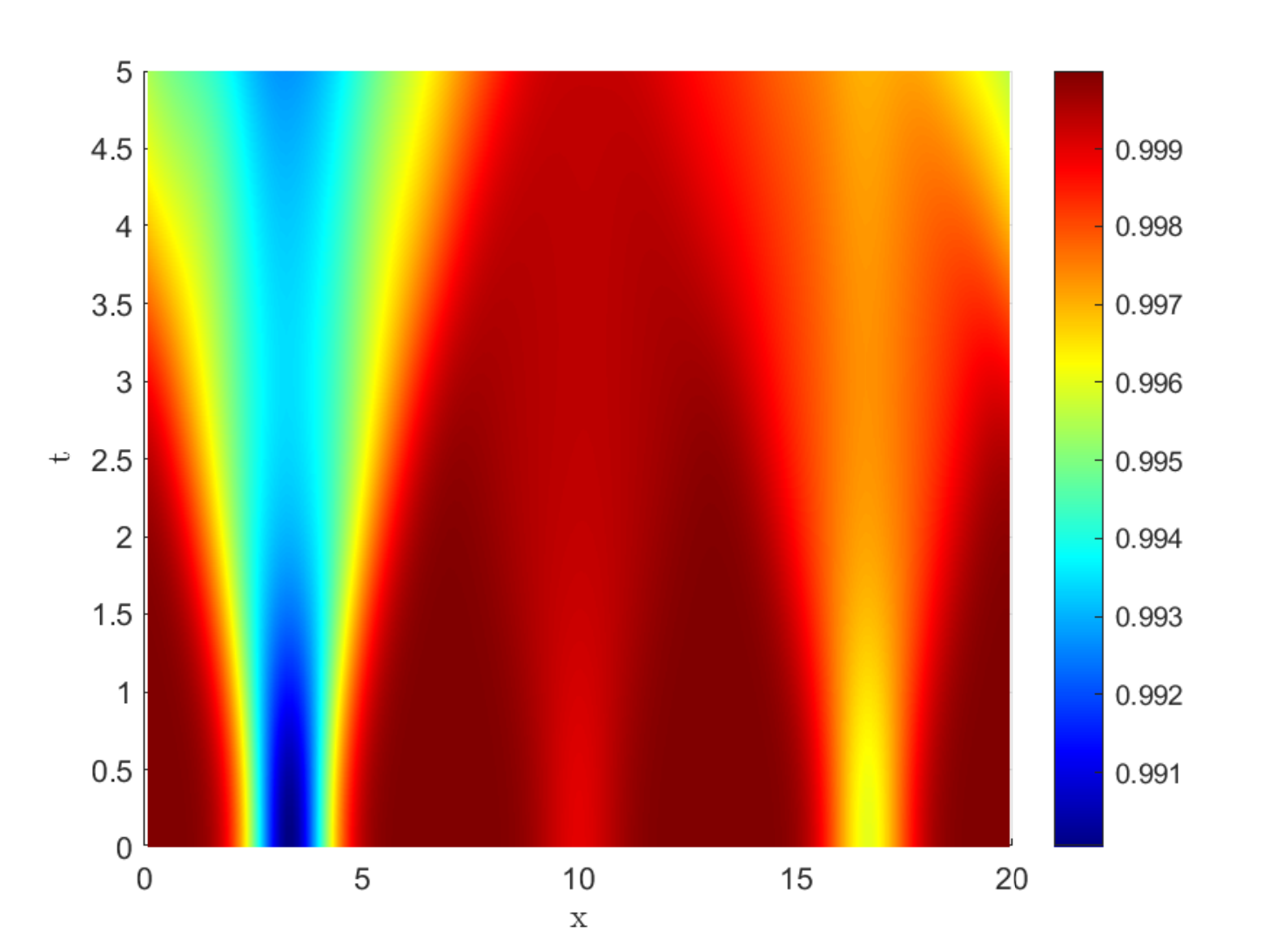}
\includegraphics[width=0.48\textwidth]{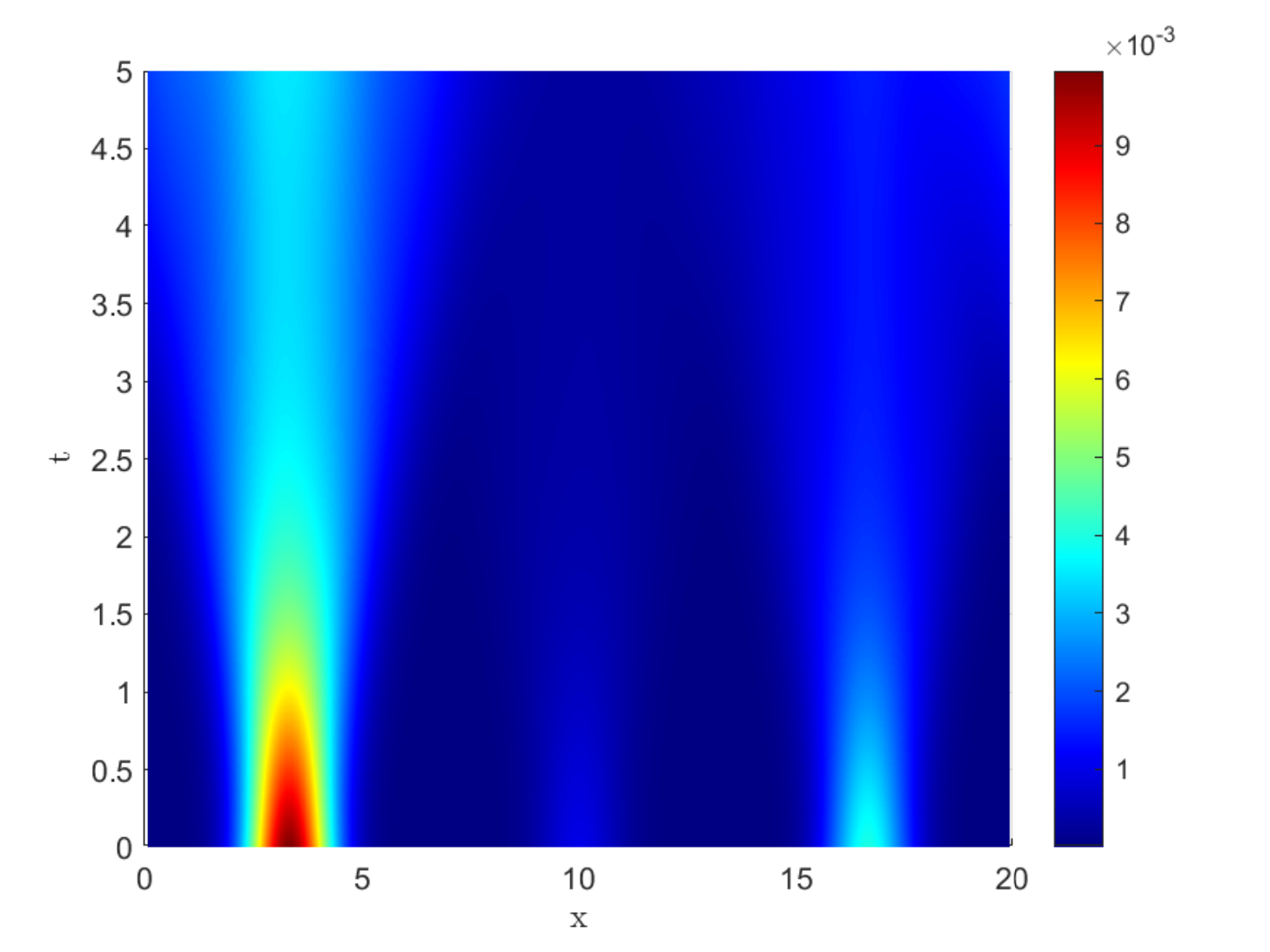}
\includegraphics[width=0.48\textwidth]{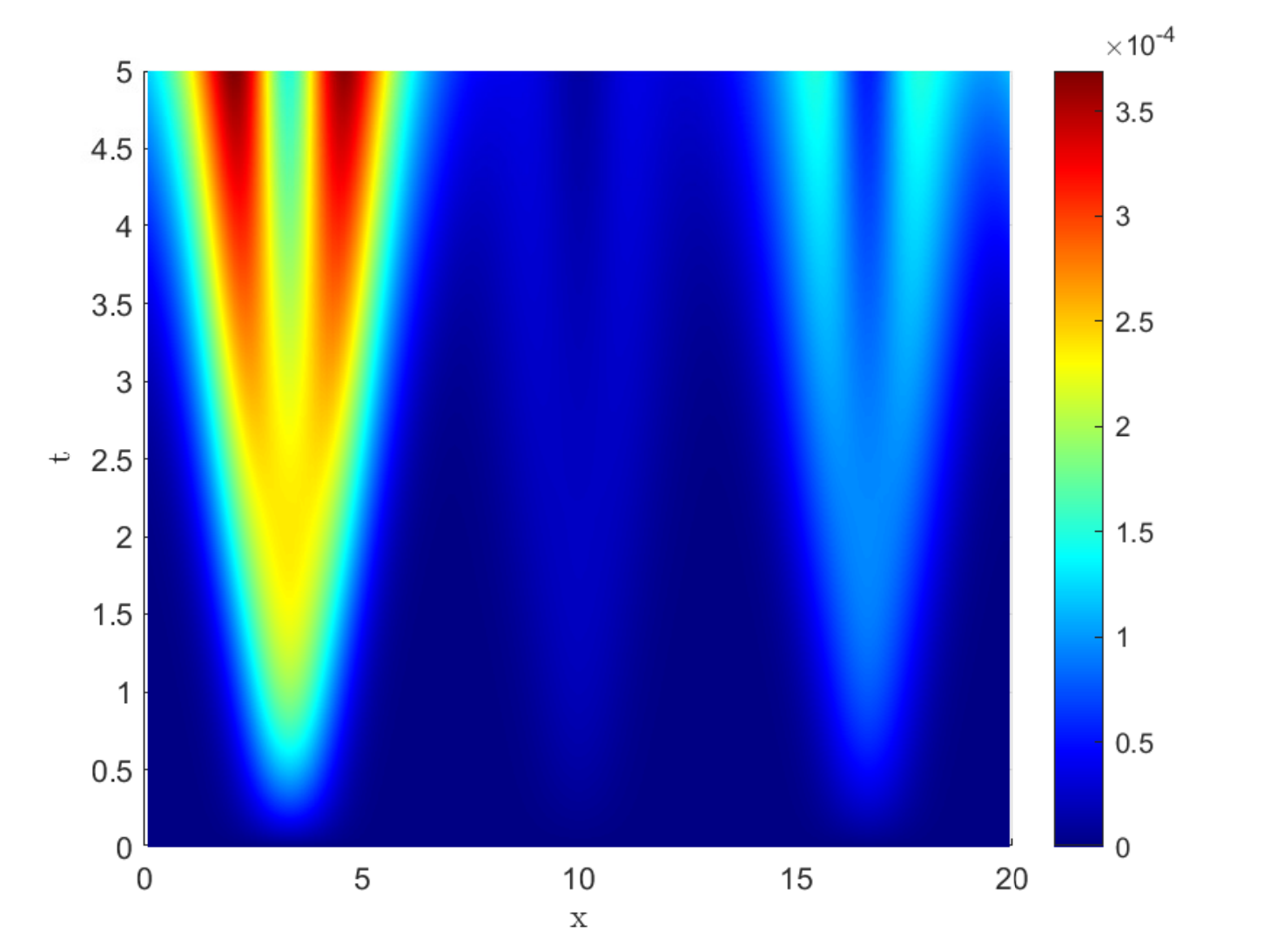}
\includegraphics[width=0.48\textwidth]{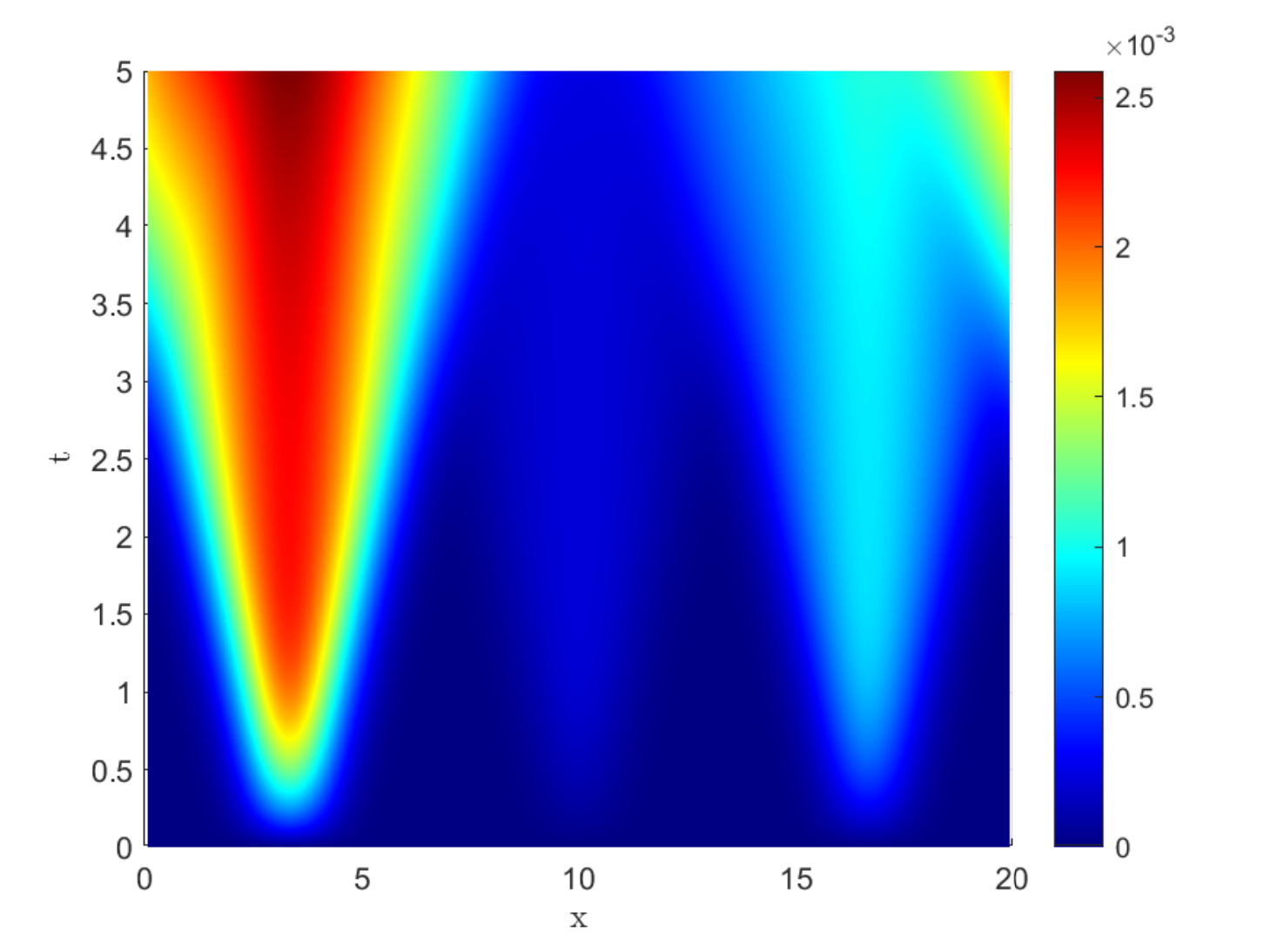}
\caption{Test 2 (b): SEIAR model in hyperbolic regime. Baseline temporal and spatial evolution of compartments $S$ (first row, left), $E$ (first row, right), $I$ (second row, left) and $A$ (second row, right) in the high-fidelity model.}
\label{fig.test2b_xt}
\end{figure}
\begin{figure}[p]
\centering
\includegraphics[scale=0.38]{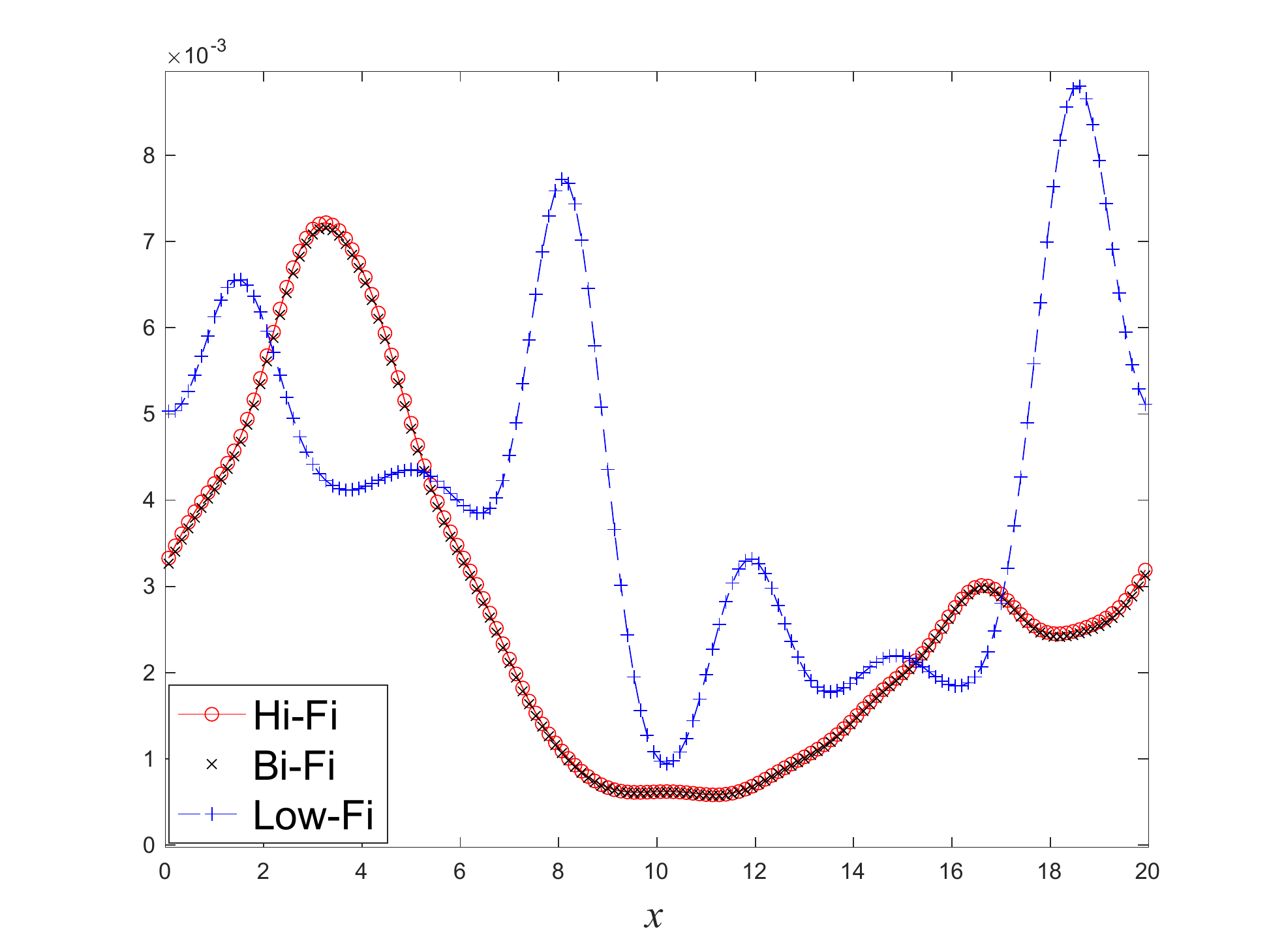}
\includegraphics[scale=0.38]{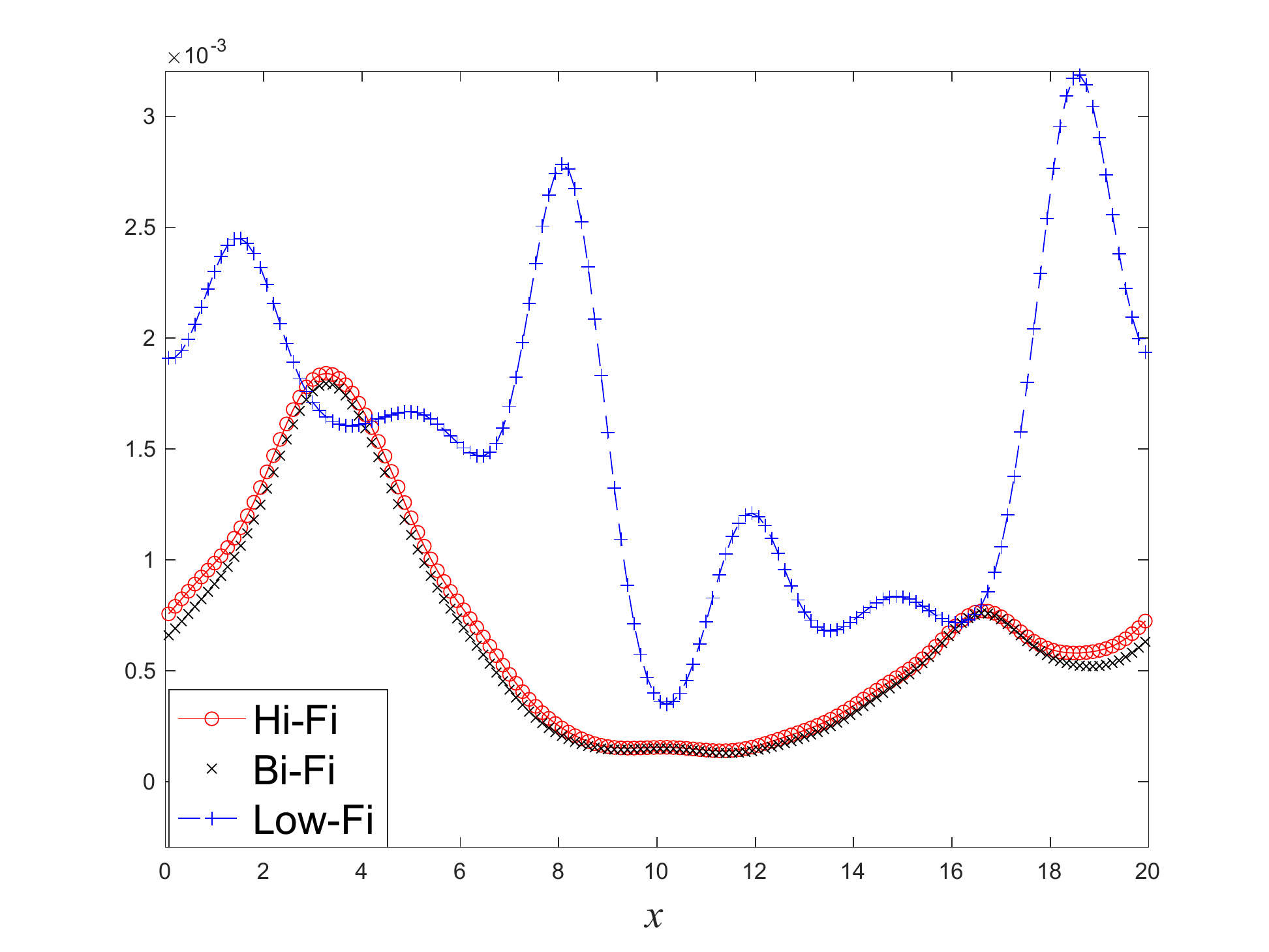}
\includegraphics[scale=0.38]{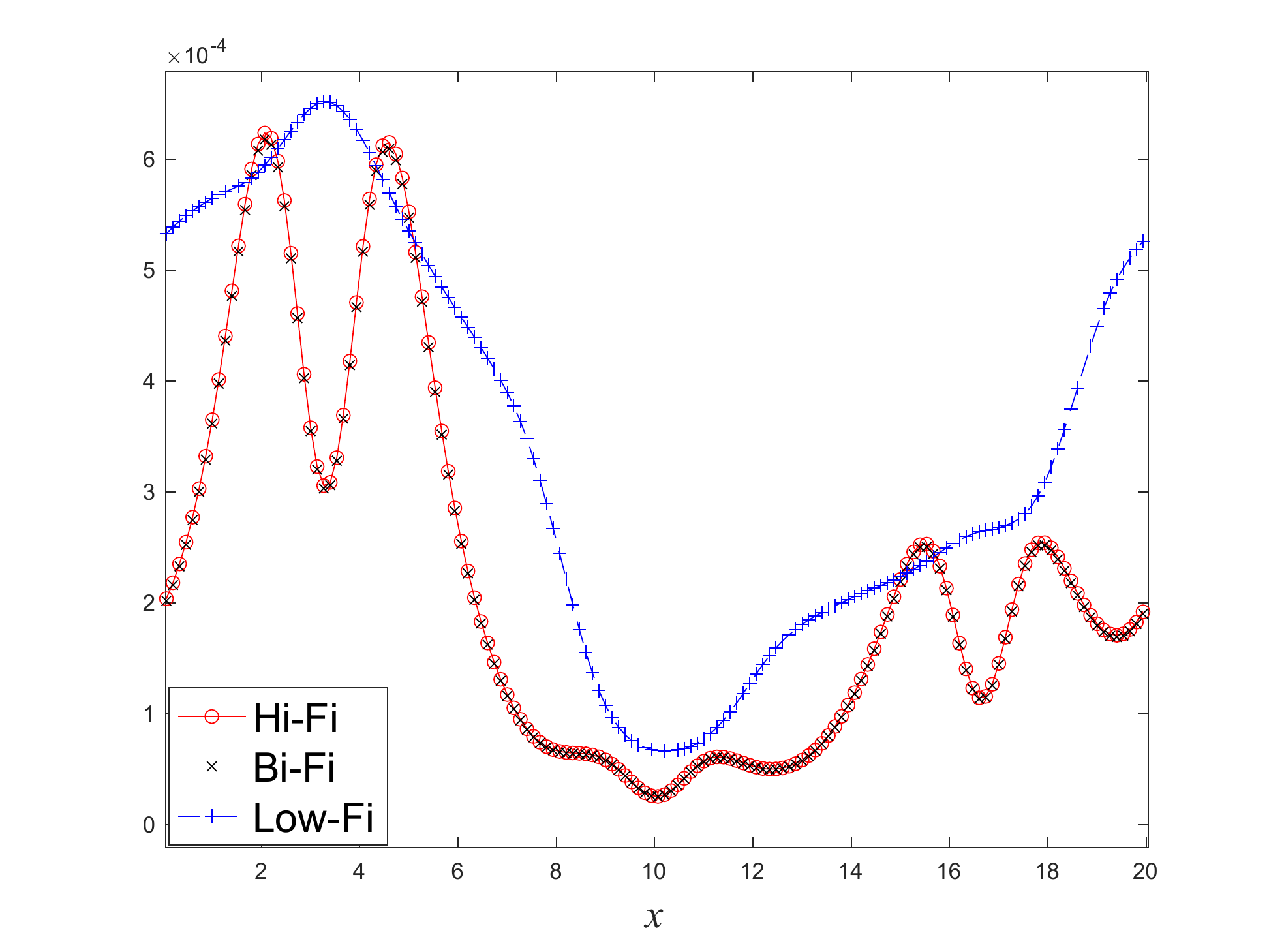}
\includegraphics[scale=0.38]{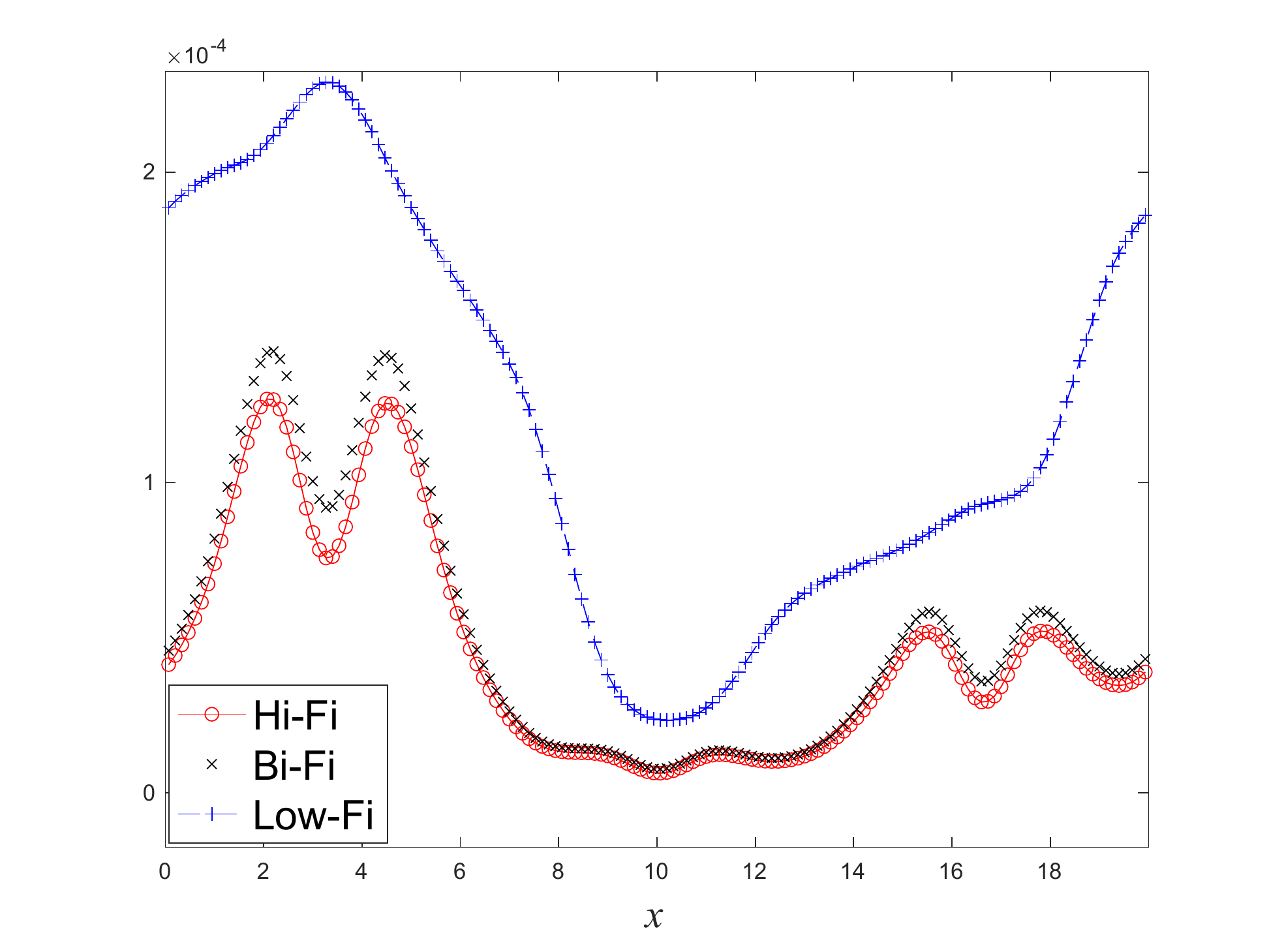}
\includegraphics[scale=0.38]{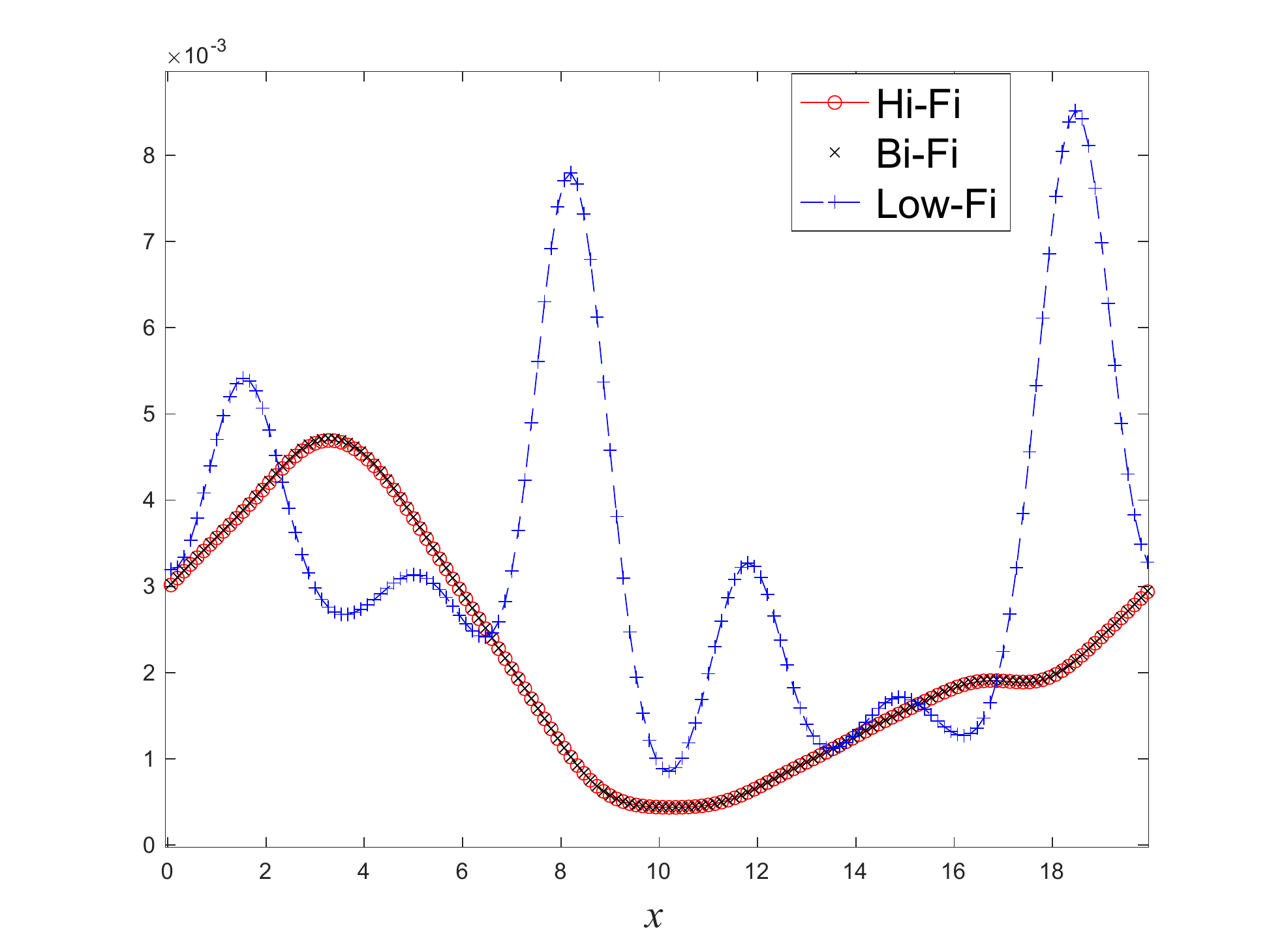}
\includegraphics[scale=0.38]{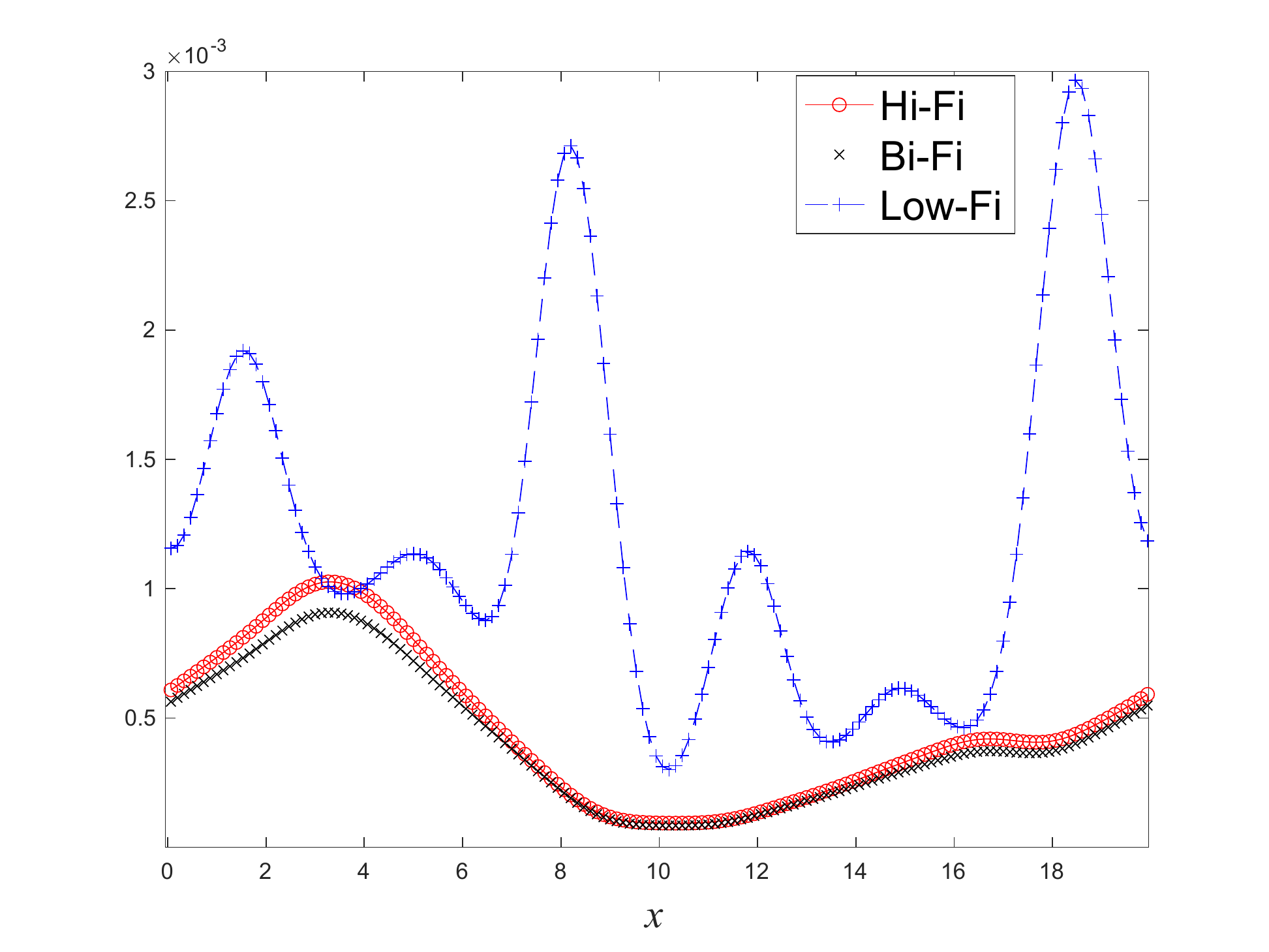}
\caption{Test 2 (b): SEIAR model in hyperbolic regime. Expectation (left) and standard deviation (right) of densities $E$ (first row), $I$ (second row) and $A$ (third row) at time $t=5$, obtained with the three methodologies, using $n=7$ for the bi-fidelity solution.}
\label{fig.test2b}
\end{figure}
\begin{figure}[t!]
\centering
\includegraphics[scale=0.33]{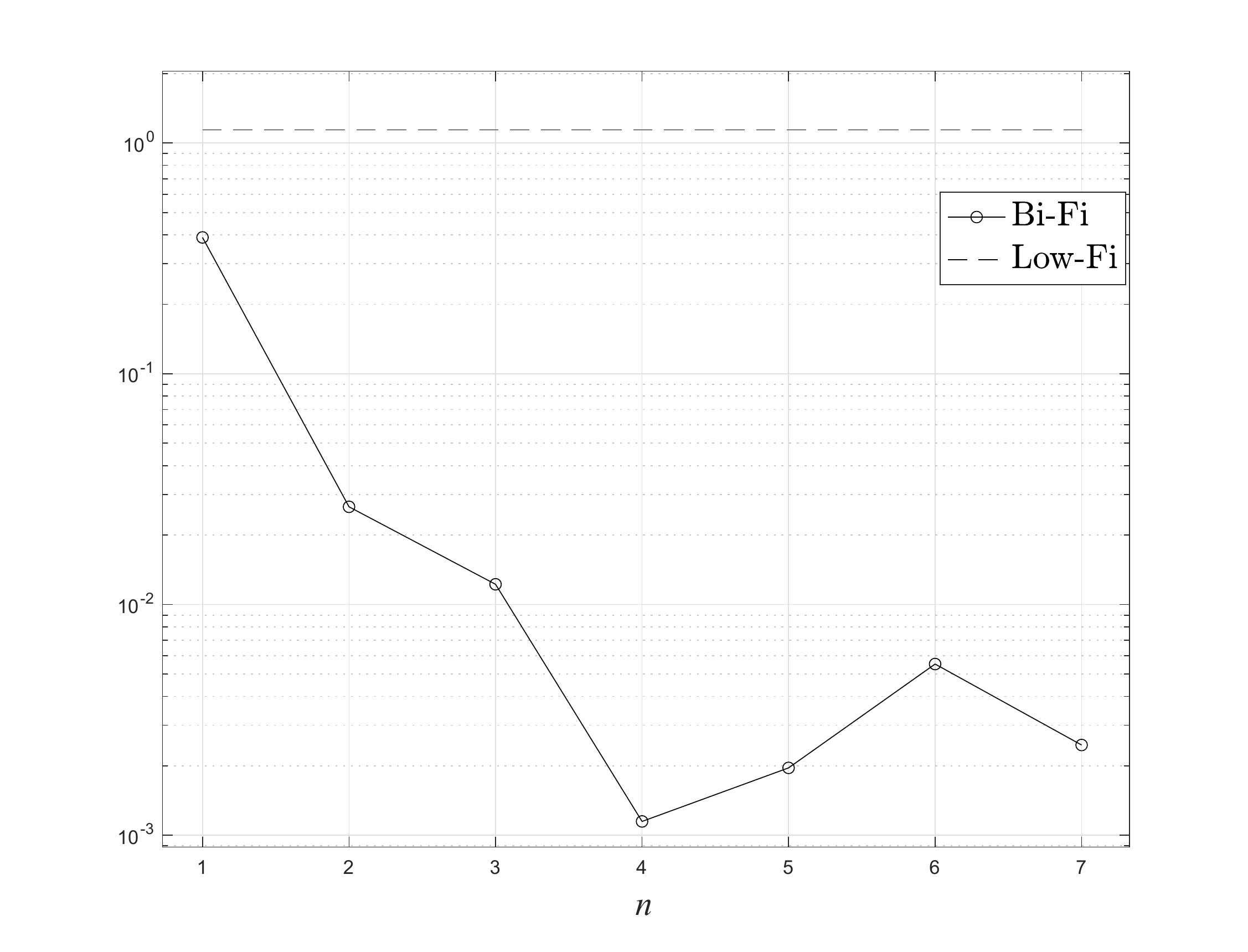}
\includegraphics[scale=0.33]{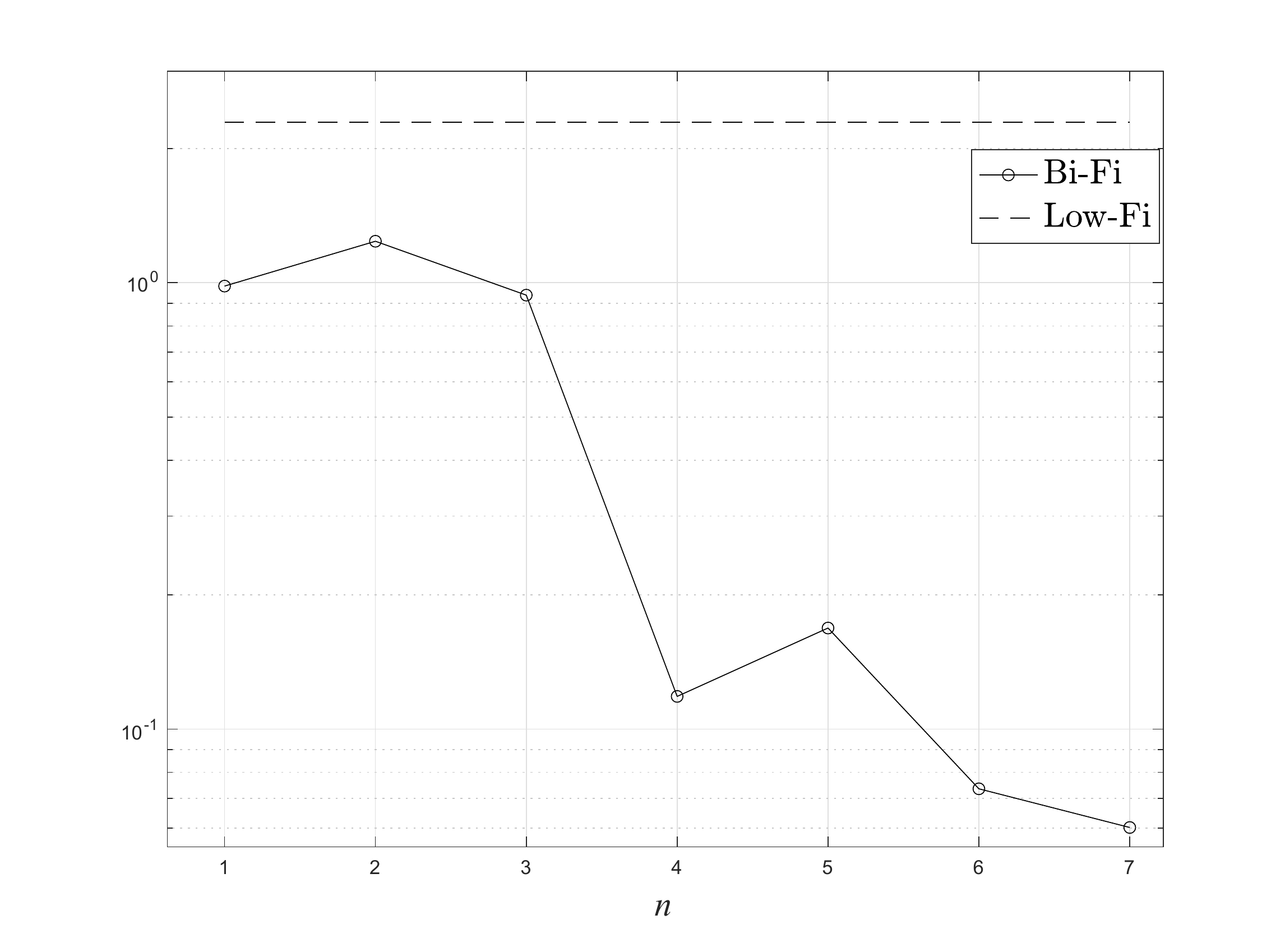}
\caption{Test 2 (b): SEIAR model in hyperbolic regime. Relative $L^2$ error decay of the bi-fidelity approximation of expectation (left) and standard deviation (right)  for the density $A$ with respect to the number of selected "important" points $n$, compared with low-fidelity errors. }
\label{fig.test2b_err}
\end{figure}

\paragraph{Test 2 (b):}
In the second case, we consider a fully hyperbolic regime choosing $\lambda_i^2=1$, $i~\in~\{S,E,A,R\}$, $\lambda_I = 0$, and $\tau_i=10$, $i~\in~\{S,E,I,A,R\}$, and $\tau_i=30$ respectively in the low-fidelity and high-fidelity model. 

The baseline space-time evolution of the compartments of major interest is presented in Figure \ref{fig.test2b_xt}. In this figure, the spatial heterogeneity of the epidemic spread related to the prescribed contact rate can be clearly appreciated. 
Mean and standard deviation obtained with the 3 approaches are shown for compartments $E$, $I$ and $A$ in Figure \ref{fig.test2b}. Again, although the low-fidelity and high-fidelity solutions show clear discordant trends, the effectiveness of the bi-fidelity method is plainly confirmed. 
In Figure \ref{fig.test2b_err}, the decay of the $L^2$ relative error norm of the bi-fidelity approximation with respect to the number of points $n$ shows that a proper accuracy can be achieved for both mean and standard deviation even in this challenging test case. 

\section{Conclusions}
In this work we introduced a bi-fidelity method for the quantification of uncertainty in epidemiological transport models based on an asymptotic-preserving space-time discretization. In detail, after presenting the high-fidelity epidemiological model, we considered the corresponding  bi-fidelity model. Both models share the same diffusive limit and permits to recover classical epidemic models based on diffusion equations. The numerical scheme used allows to have an efficient quantification of the uncertainty in different regimes, using few simulations of the high-fidelity model and several runs of the  low-fidelity model for points selection in the random space. Results for one-dimensional transport problems based on realistic compartmentalization in relation to the recent COVID-19 pandemic, which also include asymptomatic individuals, show the validity of the presented approach. Further research will be directed toward extending the present approach to realistic contexts such as those studied in^^>\cite{Bert3, Bert2}.

\section*{Acknowledgments}
G.B. and  L.P. acknowledge the support of MIUR-PRIN Project 2017, No. 2017KKJP4X \emph{Innovative numerical methods for evolutionary partial differential equations and applications}. G.B. holds a Research Fellowship from the Italian National Institute of High Mathematics, INdAM (GNCS group). L.L. holds the Direct Grant for Research supported by Chinese University of Hong Kong and Early Career Scheme 2021/22, No. 24301021, from Research Grants Council of Hong Kong. X.Z. was supported by the Simons Foundation (504054).

\end{document}